\documentclass[a4paper]{scrartcl}

\usepackage[utf8]{luainputenc}
\usepackage[USenglish]{babel}
\usepackage{csquotes}

\usepackage[a4paper, top=2.5cm, bottom=2.5cm, left=2.5cm, right=2.5cm]{geometry}
\usepackage[plainpages=false,pdfpagelabels,hidelinks,unicode]{hyperref}

\usepackage[%
  backend=bibtex,bibencoding=ascii,
  style=numeric-comp,
  giveninits=true, uniquename=init, 
  natbib=true,
  url=true,
  doi=true,
  isbn=false,
  backref=false,
  maxnames=99,
  ]{biblatex}
\usepackage{amsmath}
\numberwithin{equation}{section}
\allowdisplaybreaks
\usepackage{amssymb}
\usepackage{commath}
\usepackage{mathtools}
\usepackage{bbm}
\usepackage{nicefrac}

\usepackage{siunitx}

\usepackage{algorithm}
\usepackage{algpseudocode}

\usepackage{amsthm}
\theoremstyle{plain}
  \newtheorem{theorem}{Theorem}

\theoremstyle{definition}
  
  \newtheorem{remark}[theorem]{Remark}
  
\numberwithin{theorem}{section}

\usepackage{xcolor}
\usepackage{graphicx}
\usepackage[small]{caption}
\usepackage{subcaption}

\begingroup\expandafter\expandafter\expandafter\endgroup
\expandafter\ifx\csname pdfsuppresswarningpagegroup\endcsname\relax
\else
\fi

\usepackage{booktabs}
\usepackage{rotating}
\usepackage{multirow}

\usepackage{ifluatex}
\ifluatex
  \usepackage[no-math]{fontspec}
\else
  \usepackage[T1]{fontenc}
\fi
\usepackage{newpxtext,newpxmath}

\newcommand{\R}{\mathbb{R}}

\newcommand{\dx}{{\Delta x}}
\newcommand{\dt}{{\Delta t}}

\newcommand{\bhat}{\widehat{b}}
\newcommand{\uhat}{\widehat{u}}
\newcommand{\qhat}{\widehat{q}}
\newcommand{\tol}{\ensuremath{\texttt{tol}}}
\newcommand{\atol}{\ensuremath{\texttt{atol}}}
\newcommand{\rtol}{\ensuremath{\texttt{rtol}}}
\newcommand{\err}{\ensuremath{e}}

\newcommand{\cfl}{\nu}
\newcommand{\maxspeed}{\lambda_\mathrm{max}}

\newcommand{\fnc}[1]{\ensuremath{\mathcal{#1}}}
\newcommand{\ndof}{m}

\renewcommand{\rho}{\varrho}
\renewcommand{\epsilon}{\varepsilon}

\NewDocumentCommand{\RK}{o m O{\the\numexpr#2-1\relax} m O{} O{} o}{%
  \IfValueTF{#1}{#1}{RK}%
  #2(#3)#4%
  \ifblank{#6}{}{\textsubscript{F}}%
  \ifblank{#5}{}{[#5]}%
  \IfValueT{#7}{#7}%
}
\newcommand{\ZN}{2N}
\newcommand{\ZRp}{2R\textsubscript{+}}
\newcommand{\ERp}{3R\textsubscript{+}}
\newcommand{\ESstar}{3S*}
\newcommand{\ESstarp}{3S*\textsubscript{+}}
\newcommand{\ssp}[2]{\RK[SSP]{#2}{#1}[\ESstarp]}

\usepackage{xspace}
\newcommand{\ca}[0]{{ca.\@}\xspace}
\newcommand{\cf}[0]{{cf.\@}\xspace}
\newcommand{\eg}[0]{{e.g.\@}\xspace}
\newcommand{\ie}[0]{{i.e.\@}\xspace}

\newcommand{\mach}{\mathrm{Ma}}
\newcommand{\reynolds}{\mathrm{Re}}

\usepackage{pdfpages}

\newenvironment{keywords}{\par\textbf{Key words.}}{\par}
\newenvironment{AMS}{\par\textbf{AMS subject classification.}}{\par}

\title{Optimized Runge-Kutta Methods with Automatic Step Size Control for Compressible Computational Fluid Dynamics}
\author{Hendrik Ranocha \and
        Lisandro Dalcin \and
        Matteo Parsani \and
        David I. Ketcheson}
\date{July 21, 2021} 

\makeatletter
\hypersetup{pdfauthor={Hendrik Ranocha, Lisandro Dalcin, Matteo Parsani, David I. Ketcheson}}
\hypersetup{pdftitle={\@title}}
\makeatother

\begin{document}

\maketitle

\begin{abstract}
We develop error-control based time integration algorithms for
compressible fluid dynamics (CFD) applications and show that they are
efficient and robust in both the accuracy-limited and stability-limited
regime.  
Focusing on discontinuous spectral element semidiscretizations,
we design new controllers for existing methods and for some new embedded
Runge-Kutta pairs.
We demonstrate the importance of choosing adequate controller
parameters and provide a means to obtain these in practice.
We compare a wide range of error-control-based methods, along with the common
approach in which step size control is based on the Courant-Friedrichs-Lewy (CFL) number.  
The optimized methods give improved performance
and naturally adopt a step size close to the maximum stable CFL
number at loose tolerances, while additionally providing control of the
temporal error at tighter tolerances.  The numerical examples include
challenging industrial CFD applications.

\end{abstract}



\section{Introduction}
\label{sec:introduction}

Systems of hyperbolic conservation laws are used to model
many areas of science and engineering, such as fluid dynamics, acoustics, and
electrodynamics. In practical applications, these systems must often be solved numerically.
Explicit Runge-Kutta schemes are the most commonly used
time discretizations for hyperbolic partial differential equations (PDEs),
because of their efficiency and parallel scalability
\cite{rjl:fdmbook,kopriva2013assessment,gottlieb2016time}.
Overall efficiency of the method also depends on choosing a time step
that is as large as possible while still satisfying stability and
accuracy requirements.  Since stability requirements are frequently
more restrictive in this setting, hyperbolic PDE practitioners often adapt the
time step size based on a desired CFL number.  The CFL number involves
the ratio of the maximum characteristic speed to the mesh spacing, which is
essentially a proxy for the norm of the Jacobian.  The optimal CFL number
depends on the space and time discretizations chosen, and possibly on the problem;
it is often determined by trial and error.

On the other hand, time integration research has long emphasized the efficiency of error-based
step size control.  Much effort has gone into the design of embedded
Runge-Kutta pairs and step size controllers for this purpose.
Compared to CFL-based control, error-based control has the advantage of
not requiring a manually-tuned CFL number and allowing for control of the
temporal error when necessary.  CFL-based control has the advantage of (usually)
yielding near-optimal efficiency once the appropriate CFL value has been found,
as long as the calculation is indeed stability-limited.  
Error-based step size control for convection-dominated problems has been
attempted previously; see e.g. \cite{berzins1995temporal,ware1995adaptive}.
An ideal time integration
algorithm would achieve the efficiency of the CFL-based controller
in the stability-limited regime without the need for manually-tuned parameters,
while automatically reducing the step size if error control becomes a more
restrictive requirement.  In this work, we develop such algorithms
in the context of computational fluid dynamics (CFD).

Specifically, we focus on low-storage Runge-Kutta pairs (reviewed
in Section~\ref{sec:RK}) combined with PID step size controllers
(reviewed in Section~\ref{sec:cfl-vs-error}) and spectral element methods.
Spectral element methods can be very efficient for large-scale computations
\cite{baggag1999parallel,karniadakis2013spectral,vincent_petascale_pyfr,hutchinson2016efficiency,hadri_ccpe_2019}.
Because stability is a challenging issue for these schemes, a lot of effort has
been devoted to developing energy stable (linearly
stable) \cite{vincent_energy_fr_2011,OREILLY2017572,almquist_elastic_2020},
and entropy stable (nonlinearly stable) spatial discretizations
\cite{fisher2013high,carpenter2014entropy,gassner2016split,sjogreen2017skew,sjogreen2018high,chan2019efficient,rojas2021robustness,fernandez2020entropy}.
Stable fully-discrete schemes can be obtained from these semi-discretizations
by using a slight modification of classical time integration schemes, based on the
relaxation approach \cite{ketcheson2019relaxation,ranocha2020relaxation,ranocha2020general,ranocha2020fully}.

In the paradigm of CFL-based error control, a common approach to time integrator
design is to seek a large region of absolute stability (see \eg \cite{figueroa2021explicit}
for a recent example of this approach in the context of CFD).
For error-based control, a large region of absolute stability is again
important (for both the main method and the embedded method).  Additionally,
when automatic error control is used with step sizes near the stability
limit, the concept of step size control stability becomes crucial to
the design of the controller.  We demonstrate the importance of choosing good
step size controllers in Section~\ref{sec:importance-of-controller}.
There exists some previous work on developing error-based step size control
techniques for convection-dominated problems, principally by Berzins and co-authors \cite{berzins1995temporal,ware1995adaptive}.

We compare some existing Runge-Kutta pairs in
Section~\ref{sec:comparing-existing-schemes}, and develop optimized Runge-Kutta
pairs for discontinuous spectral element semidiscretizations of hyperbolic conservation laws
in Section~\ref{sec:optimized-RK}.
The spectral element methods applied for the numerical experiments are
implemented in the $hp$-adaptive, unstructured, curvilinear grid solver SSDC
\cite{parsani2021ssdc}.
SSDC is built on top of the Portable and Extensible Toolkit for Scientific
computing (PETSc) \cite{petsc313}, its mesh topology abstraction (DMPLEX)
\cite{knepley2009mesh}, and its scalable ODE/DAE solver library
\cite{abhyankar2018petscts}.
Further details on the spatial semidiscretizations can be found in
\cite{carpenter2014entropy,parsani2015entropy,parsani2015entropyWall,carpenter2016towards,rojas2021robustness,fernandez2020entropy}.
We perform numerical experiments
using the novel schemes in Section~\ref{sec:numerical-experiments}, both
for the compressible Euler and Navier-Stokes equations. Finally, we summarize
and discuss our results in Section~\ref{sec:summary}.
We contributed our optimized methods to the freely available open source
library DifferentialEquations.jl \cite{rackauckas2017differentialequations}
written in Julia \cite{bezanson2017julia}.

\section{Runge-Kutta methods and adaptive time stepping}
\label{sec:RK}

Using the method of lines, a spatial semidiscretization of a hyperbolic PDE
yields an ordinary differential equation (ODE) system
\begin{equation}
\label{eq:ode}
\begin{aligned}
  \od{}{t} u(t) &= f(t, u(t)),
  && t \in (0,T),
  \\
  u(0) &= u_0,
\end{aligned}
\end{equation}
where $u\colon [0,T] \to \R^\ndof$ and $\ndof$ is the number of degrees of freedom in
the spatial discretization.
An explicit first-same-as-last (FSAL) Runge-Kutta pair with $s$ stages can be
described by its Butcher tableau \cite{hairer2008solving,butcher2016numerical}
\begin{equation}
\label{eq:butcher}
\begin{array}{c | c}
  c & A
  \\ \hline
    & b^T
  \\
    & \bhat^T
\end{array}
\end{equation}
where $A \in \R^{s \times s}$ is strictly lower-triangular, $b, c \in \R^s$,
and $\bhat\in\R^{s+1}$. For \eqref{eq:ode},
a step from $u^n \approx u(t_n)$ to $u^{n+1} \approx u(t_{n+1})$, where
$t_{n+1} = t_n + \dt_n$, is given by
\begin{equation}
\label{eq:RK-step}
\begin{aligned}
  y_i
  &=
  u^n + \dt_n \sum_{j=1}^{i-1} a_{ij} \, f(t_n + c_j \dt_n, y_j),
  \qquad i \in \set{1, \dots, s},
  \\
  u^{n+1}
  &=
  u^n + \dt_n \sum_{i=1}^{s} b_{i} \, f(t_n + c_i \dt_n, y_i),
  \\
  \uhat^{n+1}
  &=
  u^n + \dt_n \sum_{i=1}^{s} \bhat_{i} \, f(t_n + c_i \dt_n, y_i) + \bhat_{s+1} f(t_{n+1}, u^{n+1}).
\end{aligned}
\end{equation}
Here, $y_i$ are the stage values of the Runge-Kutta method and the difference $u-\uhat$ is
used to estimate the local truncation error.
If $\bhat_{s+1}=0$ then \eqref{eq:RK-step} is an ordinary RK pair; otherwise it is
referred to as an FSAL RK pair.  The FSAL idea is to use the derivative of the
new solution as an additional input for the error estimator \cite{dormand1980family}.  If the step
is accepted, this costs nothing since the value $f(t_{n+1}, u^{n+1})$ must
be computed at the next step anyway.  Usually, for a main method of order $q$, the
embedded method is chosen to be of order $\qhat = q-1$;
\ie the schemes are used in local extrapolation mode.

\begin{remark}
  There are different notations for FSAL methods. A common alternative to our
  choice of using
  $A \in \R^{s \times s}$, $b, c \in \R^{s}$, and $\bhat \in \R^{s+1}$
  is to embed the baseline $s$-stage Runge-Kutta method in a method with $s+1$
  stages and Butcher coefficients
  \begin{equation}
    \tilde A = \begin{pmatrix} A & 0 \\ b^T & 0 \end{pmatrix} \in \R^{(s+1) \times (s+1)},
    \quad
    \tilde b = \begin{pmatrix} b \\ 0 \end{pmatrix} \in \R^{s+1},
    \quad
    \tilde c = \begin{pmatrix} c \\ 1 \end{pmatrix} \in \R^{s+1}.
  \end{equation}
  Then, the last row of $\tilde A$ is equal to $\tilde b$
  and $\bhat \in \R^{s+1}$ can be defined as usual.
\end{remark}

The common assumption $\sum_j a_{ij} = c_i$ is used
throughout this article. For methods with error-based step size control, the
initial step size is chosen using the algorithm described in
\cite[p.~169]{hairer2008solving}.

\subsection{Low-storage methods}
\label{sec:low-storage}

A typical RK implementation requires simultaneous storage of
all of the stages and/or their derivatives.  Each stage or derivative
occupies $m$ words; we refer to this amount of storage (sufficient for
holding a copy of the solution on the spatial grid at one point in time)
as a register.  A low-storage RK method is one that can be implemented using
only a few registers; herein we consider methods that require just three
or four registers.  Note that three registers is the fewest possible if one
requires an error estimator and the ability to reject a step.

We consider the low-storage method classes (with and without the FSAL technique):
\begin{itemize}
  \item \ESstar: three-register methods that include an error estimate
  \item \ESstarp: three-register methods that require a fourth register for the error estimate
\end{itemize}
Let $S_j$ denote a given storage register.
The class \ESstar\ methods, introduced in \cite{ketcheson2010runge}, can be implemented using only
three storage registers of size $\ndof$ if assignments of the form
\begin{equation}
  S_j \leftarrow S_j + f(t, S_j)
\end{equation}
can be made with only $\ndof + o(\ndof)$ memory. Otherwise, an additional register
is required. The \ESstar\ method family is parameterized by the coefficients
$c_i, \gamma_{1,i}, \gamma_{2,i}, \gamma_{3,i}, \beta_{i}, \delta_{i}$
and can be implemented as described in Algorithm~\ref{alg:3Sstar}.
We will also use \ESstar\ or \ESstarp\ to denote some
strong stability preserving (SSP) Runge-Kutta methods that can be
implemented using a slight modification of these algorithms, as described in
\cite{ketcheson2008highly}.

The FSAL technique has been applied to low-storage methods in
\cite{montijano2020FSAL};
these schemes append an FSAL stage to \ESstar\ methods to get more
coefficients for the embedded error estimator.

\begin{algorithm}
  \begin{algorithmic}
    \State $S_1 \gets u^n, S_2 \gets 0, S_3 \gets u^n$
    \ForAll{$i \in \set{1, \dots, s}$}
      \State $S_2 \gets S_2 + \delta_i S_1$
      \State $S_1 \gets \gamma_{1,i} S_1 + \gamma_{2,i} S_2 + \gamma_{3,i} S_3 + \beta_{i} \dt_n f(t_n + c_i \dt_n, S_1)$
    \EndFor
    \State $u^{n+1} \gets S_1$
    \State $\uhat^{n+1} \gets (S_2 + \delta_{s} S_1 + \delta_{s+1} S_3) / \bigl( \sum_{i=1}^{s+1} \delta_{i} \bigr)$
    \If{$\bhat_{s+1} \ne 0$}
      \State $\uhat^{n+1} \gets \uhat^{n+1} + \bhat_{s+1} \dt_n f(t_{n+1}, u^{n+1})$
    \EndIf
  \end{algorithmic}
  \caption{Minimum storage implementation of \ESstar\ methods.}
  \label{alg:3Sstar}
\end{algorithm}

All \ESstarp\ methods use an additional storage location for the
embedded error estimator. If the embedded method is not used, they reduce to
\ESstar\ methods without an embedded scheme. Their low-storage implementation is
delineated in Algorithm~\ref{alg:3Sstarp}.

\begin{algorithm}
  \begin{algorithmic}
    \State $S_1 \gets u^n, S_2 \gets 0, S_3 \gets u^n, S_4 \gets u^n$
    \ForAll{$i \in \set{1, \dots, s}$}
      \State $S_2 \gets S_2 + \delta_i S_1$
      \State $S_1 \gets \gamma_{1,i} S_1 + \gamma_{2,i} S_2 + \gamma_{3,i} S_3 + \beta_{i} \dt_n f(t_n + c_i \dt_n, S_1)$
      \State $S_4 \gets S_4 + \bhat_i \dt_n f(t_n + c_i \dt_n, S_1)$
    \EndFor
    \State $u^{n+1} \gets S_1$
    \State $\uhat^{n+1} \gets S_4$
    \If{$\bhat_{s+1} \ne 0$}
      \State $\uhat^{n+1} \gets \uhat^{n+1} + \bhat_{s+1} \dt_n f(t_{n+1}, u^{n+1})$
    \EndIf
  \end{algorithmic}
  \caption{Minimum storage implementation of \ESstarp\ methods.}
  \label{alg:3Sstarp}
\end{algorithm}

\subsection{Error-based step size control}
\label{subsec:err_step_size}

We use step size controllers based on digital signal processing
\cite{gustafsson1988pi,gustafsson1991control,soderlind2002automatic,soderlind2003digital,soderlind2006time} implemented in PETSc
\cite{petsc313,abhyankar2018petscts}. In particular, we use PID controllers
that select a new time step using the formula\footnote{Open source libraries
such as PETSc \cite{abhyankar2018petscts} and DifferentialEquations.jl,
\cite{rackauckas2017differentialequations} where we implemented a PID step size
control, will usually limit the factor multiplying the time step size, \eg using
a limiter of the form $f(x) = 1 + \arctan(x - 1)$ \cite{soderlind2006adaptive}.}
\begin{equation}
\label{eq:PID}
  \dt_{n+1} = \epsilon_{n+1}^{\beta_1 / k}
              \epsilon_{n }^{\beta_2 / k}
              \epsilon_{n-1}^{\beta_3 / k} \dt_{n},
\end{equation}
where $q$ is the order of the main method, $\qhat$ is the order of the
embedded method (usually $\qhat = q  - 1$), $k = \min(q, \qhat) + 1$
(usually $k = \qhat + 1 = q$), $\beta_i$ are the controller parameters, and
\begin{equation}
  \epsilon_{n+1} = \frac{1}{w_{n+1}},
  \quad
  w_{n+1} = \left( \frac{1}{\ndof} \sum_{i=1}^{\ndof} \left( \frac{u_i^{n+1} - \uhat_i^{n+1}}{\atol + \rtol \max\{ |u_i^{n+1}|, |\uhat_i^{n+1}| \}} \right)^2 \right)^{1/2},
\end{equation}
where $\ndof$ is the number of degrees of freedom in $u$, and $\atol$, $\rtol$
are the absolute and relative error tolerances.
Some common controller parameters recommended in the literature are given in
Table~\ref{tab:classical-controllers}.
Unless stated otherwise, we use equal absolute and relative error tolerances.
The choice of the weighted/relative error estimate $w_{n+1}$ is common in the
literature \cite[Equations (4.10) and (4.11)]{hairer2008solving} and often
the default choice in general purpose ODE software such as PETSc
\cite{abhyankar2018petscts} or DifferentialEquations.jl
\cite{rackauckas2017differentialequations}. This choice of $w_{n+1}$ allows to
decouple the time integration parameters from a possible spatial semidiscretization.
In contrast to a quadrature-based approach, it weighs degrees of freedom of
different refinement levels in the same way, which can be beneficial, since
refined regions (of interest) are not weighed less than coarse regions (without
interesting solution features).

If the factor multiplying the old time step $\dt_{n}$ is too small or the
solution is out of physical bounds, \eg because of negative density/pressure
in CFD, the step is rejected and retried with a smaller time step $\dt_{n}$.
The default options used in all numerical experiments described in this work
accept a step if the factor multiplying the step size is at least $0.9^2$.
Otherwise, the step is rejected and retried with the step size predicted by the
PID controller. If the solution is out of physical bounds, the step is rejected
and retried with a time step reduced by a factor of four.

\begin{table}[htb]
\centering
\caption{Classical step size controllers recommended in the literature.}
\label{tab:classical-controllers}
\begin{tabular}{*4c}
  \toprule
  Controller & $\beta_1$ & $\beta_2$ & $\beta_3$ \\
  \midrule
  PI42 & \num{0.60} & \num{-0.20} & \num{0.00} \\
  PI33 & \num{0.66} & \num{-0.33} & \num{0.00} \\
  PI34 & \num{0.70} & \num{-0.40} & \num{0.00} \\
  \bottomrule
\end{tabular}
\end{table}

\section{CFL- vs. error-based step size control}
\label{sec:cfl-vs-error}

The error-based step size control described above is efficient if
the practical time step is limited by the constraint of accuracy.
On the other hand, if the allowable time step is determined by stability,
and an explicit time discretization is employed, then it is natural to
use a step size of the form $\dt_n \propto 1/L_n$, where $L_n$ is
an approximation of the norm of the Jacobian of the ODE system.

In the time integration of hyperbolic PDEs, it is indeed often the case that
the step size is limited in practice by stability rather than accuracy.
Therefore it is common practice to use a step size control of the kind
just described.  For such systems, the norm of the Jacobian is proportional
to $\max_i(\maxspeed(u^n_i)/\dx_i)$, where $\dx_i$ is a local measure of the mesh
spacing (at grid point/cell/element $i$), and
$\maxspeed$ is the maximal (local) wave speed, related to the largest-magnitude
eigenvalue of the flux Jacobian of the hyperbolic system.
The step size control thus takes the form (referred to herein as a CFL-based control)
\begin{equation}
\label{eq:cfl-dt}
  \dt_n = \cfl \, \min_i \frac{\dx_i}{\maxspeed(u^n_i)},
\end{equation}
where $\cfl$ is the desired CFL number.  The appropriate choice of $\cfl$
depends on the details of the space and time discretizations; it can be
studied theoretically using linearization (see \eg \cite{langseth2000wave})
but is often determined experimentally.  
An additional complication
is the question of how to define $\dx$. Even on uniform Cartesian grids and regular triangulations
\cite{kubatko2008time}, multiple waves traveling in different directions make an optimal choice of $\cfl$ difficult.
This is even more a challenging question for unstructured grids.

A clear advantage of error-based control is the availability of an
estimate of the temporal error.
At first glance, error-based step size control seems inappropriate in
the stability-limited regime, since the local error may not be very sensitive to
small differences between stable and unstable step sizes, near the stable
step size limit.  A tight error tolerance that ensures stability at all steps
might result in an excessively small step size.
However, as described in \cite[Section~IV.2]{hairer2010solving} and
discussed below, it is possible to design error-based step size controllers
that behave appropriately in the stability-limited regime.

Both classes of controllers require some user-determined parameters:
$\cfl$ and $\dx$ for the CFL-based controller, and $\atol$ and $\rtol$ for the error-based
controller.  In this section we show through an example that
carefully-designed error-based controllers can achieve
near-maximal efficiency in a way that is relatively insensitive to changes
in the user parameters.  In contrast, the efficiency of the CFL-based controller
always bears a linear sensitivity to the parameters $\cfl$ and $\dx$.

To demonstrate, we consider the two-dimensional
advection equation with constant velocity $a = (1, 1)^T$
in the domain $[-5,5]^2$ with periodic boundary conditions.
An initial sinusoid of one wavelength in each direction
is advected over the time interval $[0, 100]$.
In space we apply the spectral collocation method of SSDC based on solution polynomials
of degree $p = 4$ \cite{parsani2021ssdc}.

For the CFL-based controller, the ratio of the local mesh spacing and the
maximal speed at a node $i$ is estimated in this case as
\begin{equation}
\label{eq:cfl-factor}
  \frac{\dx_i}{\maxspeed(u^n_i)}
  =
  \sigma \frac{J_i}{\sum_{j=1}^d | (J \partial_{x} \xi^j)_i \cdot a |},
\end{equation}
where $d$ is the number of spatial dimensions ($d = 2$ for this example),
$a = (1, 1)^T$ is the constant advection velocity, $J_i$ the determinant of
the grid Jacobian $\partial_{x} \xi$ at node $i$, $(J \partial_{x} \xi^j)_i$
is the contravariant basis vector in direction $j$ at node $i$
\cite[Chapter~6]{kopriva2009implementing}, and $\sigma$ is a
normalizing factor depending on the solution polynomial degree, $p$, which is usually chosen
such that a real stability interval of $2$ corresponds to $\cfl = 1$.

\begin{figure}[htb]
\centering
  \includegraphics[width=0.3\textwidth]{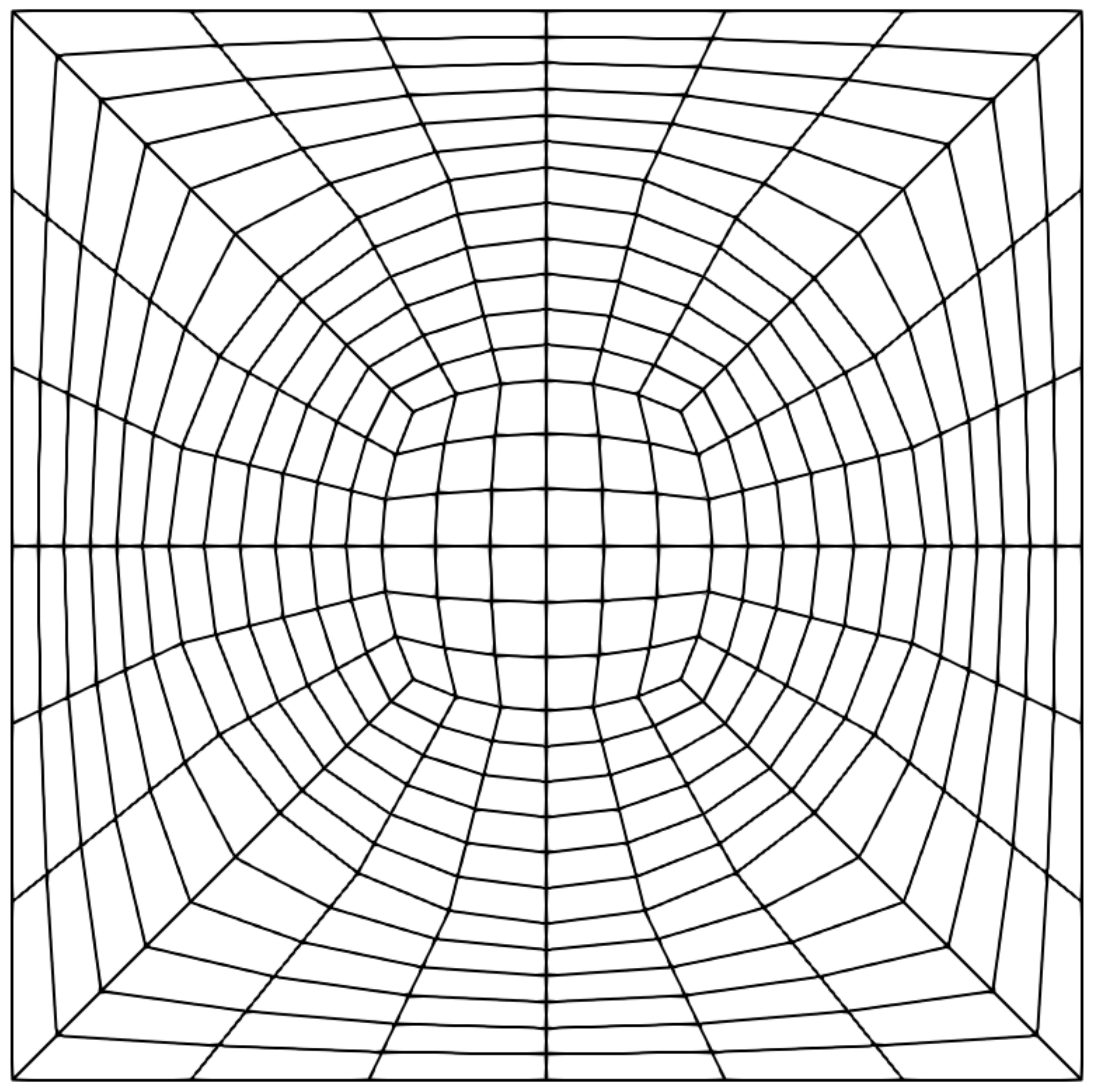}
  \caption{Unstructured grid used in the comparison of CFL- and error-based
           step size control.}
  \label{fig:mesh-box-2d}
\end{figure}

\begin{figure}[!htb]
\centering
  \begin{subfigure}{\textwidth}
  \centering
    \includegraphics[width=\textwidth]{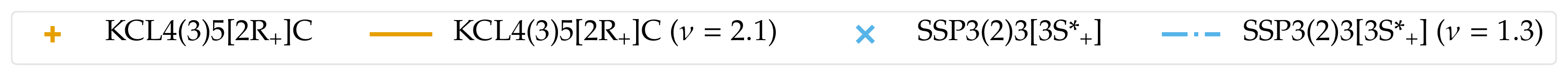}
  \end{subfigure}%
  \\
  \begin{subfigure}{0.49\textwidth}
  \centering
    \includegraphics[width=\textwidth]{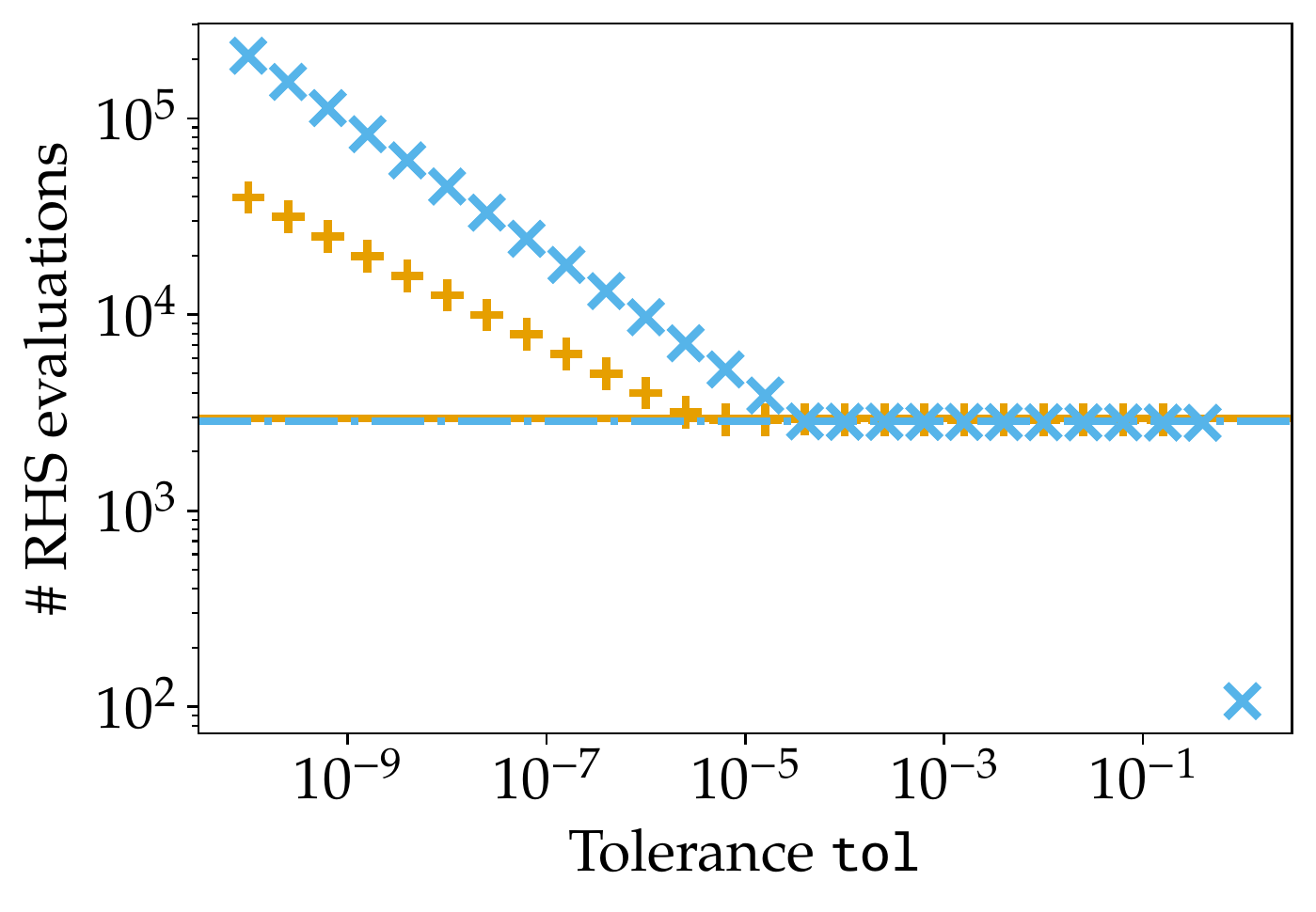}
    \caption{Number of RHS evaluations.}
  \end{subfigure}%
  \hspace*{\fill}
  \begin{subfigure}{0.49\textwidth}
  \centering
    \includegraphics[width=\textwidth]{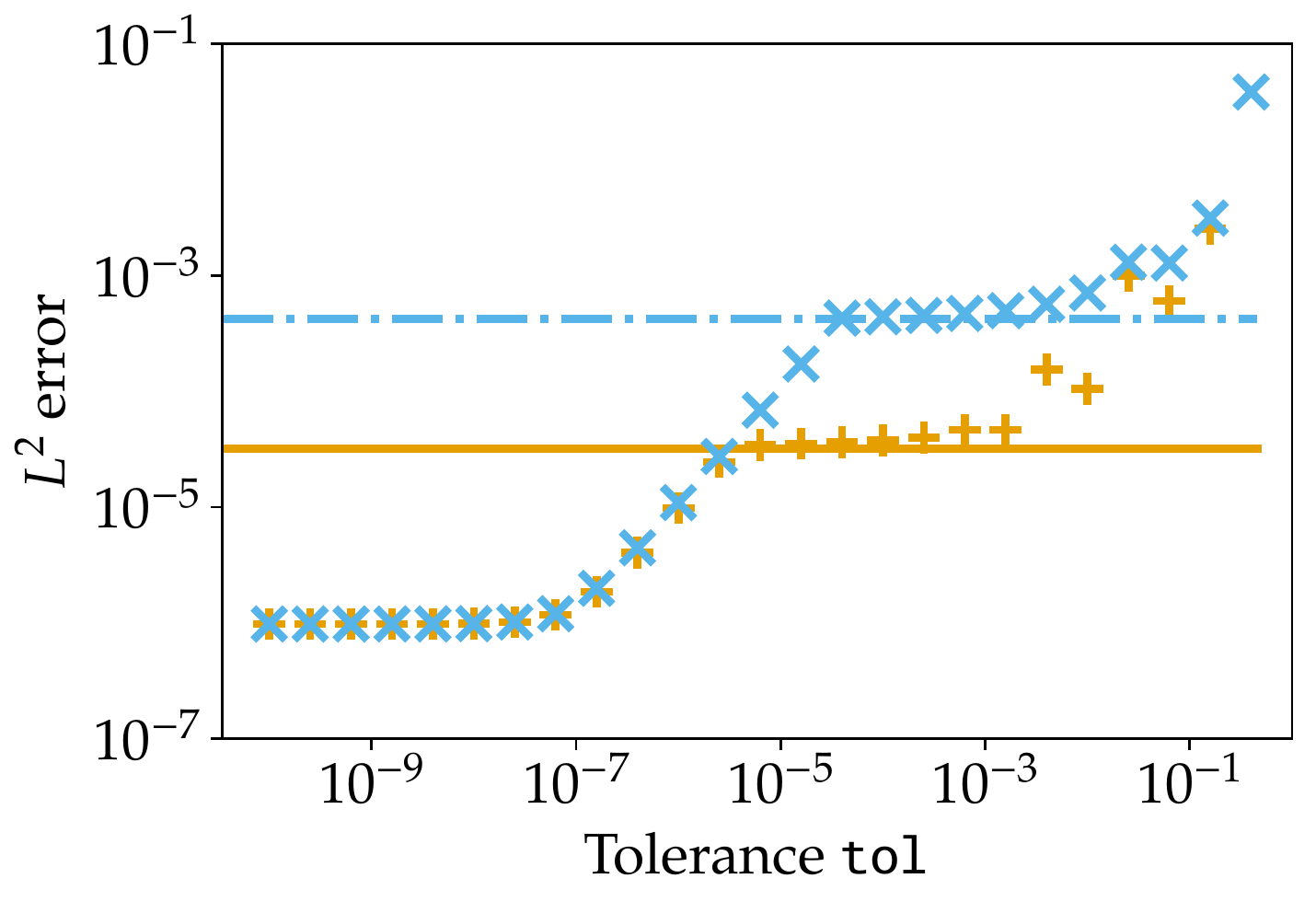}
    \caption{$L^2$ error at the final time.}
  \end{subfigure}%
  \caption{Performance of CFL- and error-based step size control for a linear
           advection problem with $p = 4$ on a uniform grid. For CFL-based
           controllers, the maximal CFL number is always included.
           The error-based controller is a standard PI controller with
           $\beta_1 = 0.7, \beta_2 = -0.4$ and uses equal absolute and
           relative tolerances $\atol = \rtol = \tol$.}
  \label{fig:CFL-vs-error-uniform}
\end{figure}

\begin{figure}[!htb]
\centering
  \begin{subfigure}{0.85\textwidth}
  \centering
    \includegraphics[width=\textwidth]{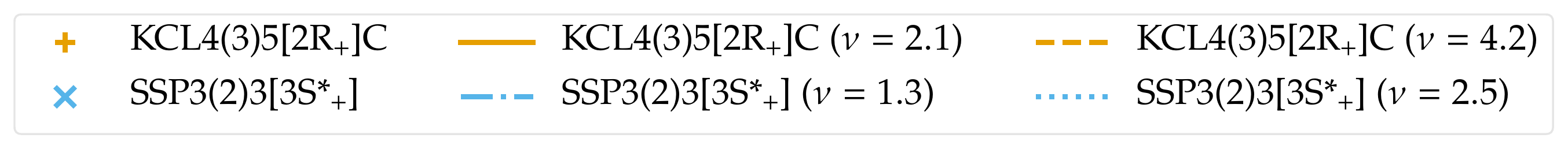}
  \end{subfigure}%
  \\
  \begin{subfigure}{0.49\textwidth}
  \centering
    \includegraphics[width=\textwidth]{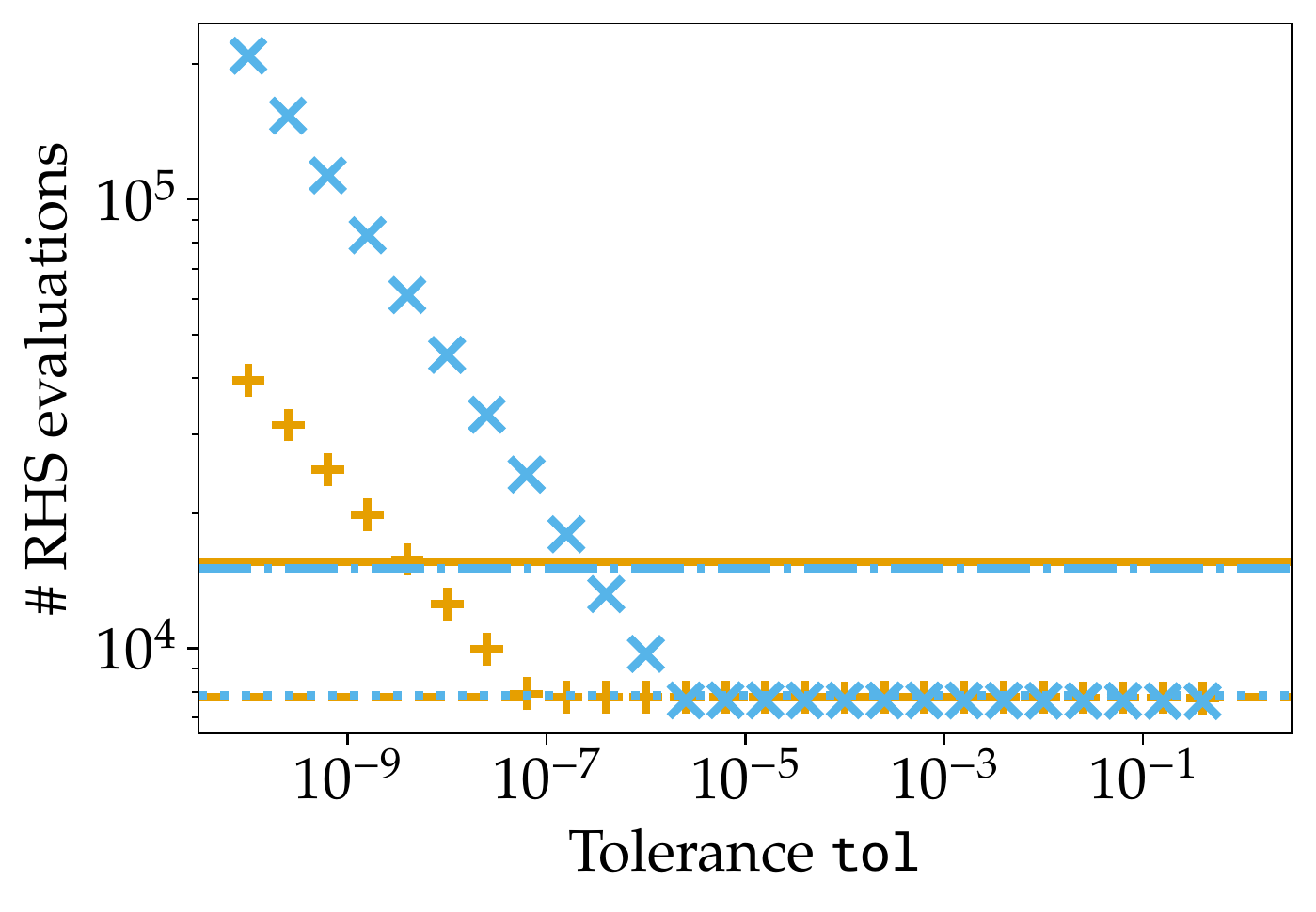}
    \caption{Number of RHS evaluations.}
  \end{subfigure}%
  \hspace*{\fill}
  \begin{subfigure}{0.49\textwidth}
  \centering
    \includegraphics[width=\textwidth]{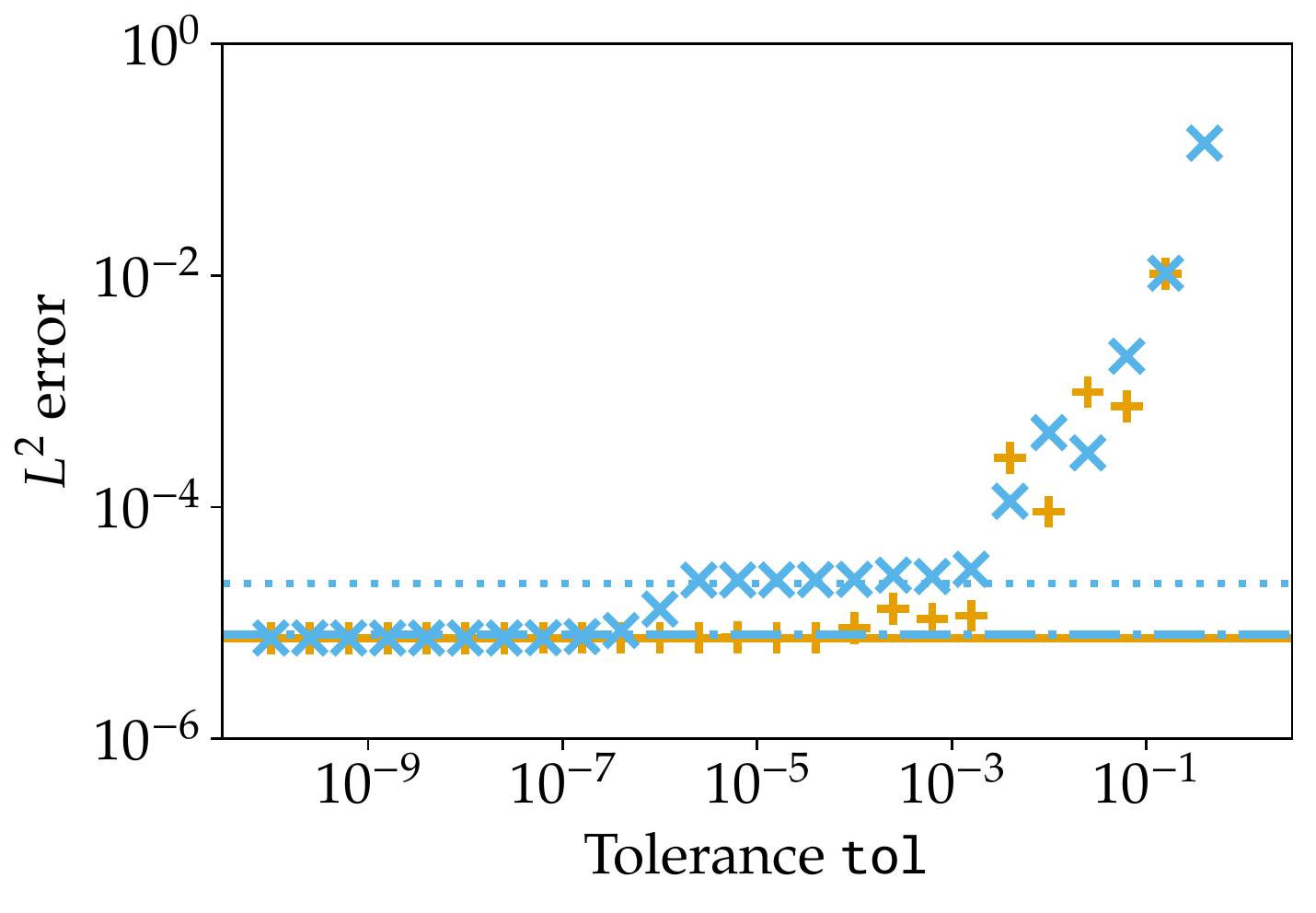}
    \caption{$L^2$ error at the final time.}
  \end{subfigure}%
  \caption{Performance of CFL- and error-based step size control for a linear
           advection problem with $p = 4$ on a non-uniform grid. For CFL-based
           controllers, the maximal CFL number is always included.
           The error-based controller is a standard PI controller with
           $\beta_1 = 0.7, \beta_2 = -0.4$ and uses equal absolute and
           relative tolerances $\atol = \rtol = \tol$.}
  \label{fig:CFL-vs-error-nonuniform}
\end{figure}

Two grids will be used: A regular, uniform grid with $8^2$ elements and the
curvilinear grid shown in Figure~\ref{fig:mesh-box-2d}.
The results are summarized in Figure~\ref{fig:CFL-vs-error-uniform}
for the uniform grid and in Figure~\ref{fig:CFL-vs-error-nonuniform}
for the unstructured grid.
We test three time discretizations: the popular fourth-order, five-stage,
low-storage method \RK[CK]{4}{5}[\ZN] of \cite{carpenter1994fourth};
the method \RK[KCL]{4}{5}[\ZRp][][C] of \cite{kennedy2000low}, which comes with an embedded third-order error estimator;
and the strong stability preserving (SSP) method \ssp33 of \cite{shu1988efficient}, equipped with the embedded method of \cite{conde2018embedded}.

\begin{remark}
  We use the same naming convention as \cite{kennedy2000low}, referring to an
  $s$-stage Runge-Kutta method of order $q$ with embedded method of order $\qhat$
  as \RK[NAME]{$q$}[$\qhat$]{$s$}.
  Additional identifiers indicating low-storage requirements or other properties
  are appended, \eg a subscript ``F'' for FSAL methods.
  The number of stages $s$ denotes the effective number of RHS evaluations
  per step, which is one less than the number of stages for FSAL methods.
  For low-storage methods, the required amount of memory based on certain
  assumptions is listed following the notation of \cite{ketcheson2010runge}.
  In particular, $n$N methods need only $n$ memory registers of size $m$ if
  assignments of the form $S_j \leftarrow \alpha S_j + f(t, S_i)$ can be made
  without additional allocations. Similarly, $nR$ methods use $n$ memory registers
  and assignments of the form $S_j \leftarrow f(t, S_j)$; $mS$ methods were
  described in Section~\ref{sec:low-storage}. As described there, a subscript
  $_+$ indicates methods that require an additional storage register if an
  embedded error estimator is used.
  Additional parts of the names of RK methods are usually taken directly from
  their sources. For example, the method \RK[KCL]{4}{5}[\ZRp][][C] is a
  fourth-order method with embedded third-order error estimator. It has five
  stages and requires two memory registers based on the $n$R assumption. If
  the embedded error estimator is used, it requires three memory registers.
  The C suffix is appended as suggested in \cite{kennedy2000low} to indicate
  a particular design criterion (in this particular case, looking for a compromise
  between linear stability and accuracy).
\end{remark}

The widely-used method \RK[CK]{4}{5}[\ZN] has linear stability properties very
similar to \RK[KCL]{4}{5}[\ZRp][][C] --- both methods have the same maximum
stable CFL number $\cfl = 2.1$. Thus, they use the same number of RHS evaluations
while yielding nearly the same errors. Thus, we omit the method \RK[CK]{4}{5}[\ZN]
in the following plots and use only the \RK[KCL]{4}{5}[\ZRp][][C] pair, for which
it is possible to use error-based step size control.
Impressively, the error-based controller manages, for a wide range of tolerances,
to use almost exactly the same number of steps as the carefully tuned CFL-based
controller.  Over this range of tolerances, including $\tol \in [10^{-5}, 10^{-3}]$,
the step size is determined by stability. Hence, the number of RHS evaluations
and the error are nearly independent of the tolerance in this regime.
For tolerances larger than $10^{-3}$, the final error increases for this
long-time simulation while the number of RHS evaluations stays nearly
the same. For tighter tolerances (below $10^{-6}$), the error-based controller
detects accuracy restrictions and increases the number of RHS evaluations.
This also leads to a reduction of the final error until it plateaus again
because of the dominant spatial error (at \ca $\tol = 10^{-7}$). However,
the number of RHS evaluations keeps increasing.

Using the same CFL number of $\cfl = 2.1$ on the unstructured grid still results
in a stable simulation. However, the CFL number can be doubled there without
increasing the error significantly.  Hence, the user has to tune this parameter
carefully to get a stable and efficient simulation.
In contrast, using error-based step size control we see behavior very similar
to what was observed for the uniform grid.
The same error tolerance can be used, resulting in the same optimal number of function evaluations
determined manually for the CFL-based step size controller. This demonstrates
the enhanced robustness properties of error-based step size control.

These examples suggest that error-based control is more robust to changes in
the grid and less sensitive to the required user parameters. Similar
results have been obtained using other Runge-Kutta schemes for this problem
and for more challenging problems, some of which are presented later in this work.
For practitioners whose primary interest is in
applying the schemes to solve challenging scientific problems or developing
spatial semidiscretizations, error-based
time step controllers seem favorable, since the most important design choices
have to be provided by the developers of the time integration schemes and
the practitioners have to choose only the rather robust error tolerance of
the solver.

\section{Importance of controller parameters}
\label{sec:importance-of-controller}

Standard error-based controllers will often work acceptably in the asymptotic
regime (\ie, the regime where the leading truncation error term strongly dominates
all subsequent terms). However, as demonstrated in Section~\ref{sec:cfl-vs-error},
applications involving convection-dominated problems are often constrained by
stability, so that one may be working outside the asymptotic regime.
In this case, the standard theory does not apply; instead, step size control stability
has to be considered \cite{hall1988analysis}.

Following \cite[Section~IV.2]{hairer2010solving},  step size control stability
can be explained using the linear model problem $\od{}{t} u(t) = \lambda u(t)$.
Given an explicit Runge-Kutta method with embedded error estimator and a
PID controller \eqref{eq:PID} with parameters $\beta_i$, the update formulae become
\begin{equation}
\begin{aligned}
  u^{n+1} &= R(\dt_n \lambda) u^n,
  \\
  \err^{n+1} &= E(\dt_n \lambda) u^n,
  \\
  \dt_{n+1} &= \text{PID}(\beta, \atol, \rtol, \err^{n+1}, \err^{n}, \err^{n-1}),
\end{aligned}
\end{equation}
where $R$ is the stability polynomial of the main method, $E$ is the difference
of the stability polynomials of the embedded and the main method, and $\err$
is the (local) error estimate.
By taking logarithms, this update formula can be reduced to a difference
recursion with fixed points on the boundary of the stability region of the
main method. To get a stable behavior, the spectral radius of the associated
Jacobian has to be less than unity \cite[Proposition~IV.2.3]{hairer2010solving}.
For a PID controller \eqref{eq:PID}, this Jacobian becomes \cite{kennedy2000low}
\begin{equation}
\label{eq:jacobian-control-stability}
\renewcommand\arraystretch{1.5}
  J(z)
  =
  \begin{psmallmatrix}
    1 & \frac{z R'(z)}{R(z)} & 0 & 0 & 0 & 0 \\
    -\frac{\beta_1}{k} & 1 - \frac{\beta_1}{k} \frac{z E'(z)}{E(z)} &
      -\frac{\beta_2}{k} & -\frac{\beta_2}{k} \frac{z E'(z)}{E(z)} &
      -\frac{\beta_3}{k} & -\frac{\beta_3}{k} \frac{z E'(z)}{E(z)} \\
    1 & 0 & 0 & 0 & 0 & 0 \\
    0 & 1 & 0 & 0 & 0 & 0 \\
    0 & 0 & 1 & 0 & 0 & 0 \\
    0 & 0 & 0 & 1 & 0 & 0 \\
  \end{psmallmatrix},
\end{equation}
where $k$ is the order of the error estimator; if $\qhat = q - 1$ is the order
of the embedded method, then $k = q = \qhat + 1$.
To get step size control stability, one can fix a controller such as the standard
I~controller and optimize the RK pair accordingly as demonstrated in
\cite{higham1990embedded}. The other possibility, pursued here, is to optimize
the controller parameters for a given RK pair.

While one might hope that a controller designed
to work well with one method will also work well with other methods,
this is generally not the case.
Rather, a controller should be designed for the given error
estimator; \cf \cite{arevalo2020local} for the case of linear multistep methods.
To demonstrate this, we consider again the linear advection problem described
in Section~\ref{sec:cfl-vs-error} with a uniform mesh.
We will take the PI34 controller with $\beta_1 = 0.7, \beta_2 = -0.4$ \cite{gustafsson1991control},
designed for use with the classical \RK[DP]{5}{6}[][FSAL] method of \cite{prince1981high},
but use instead the \RK[BS]{5}{7}[][FSAL] method of \cite{bogacki1996efficient}.  Note that
both are 5(4) pairs designed with similar purposes in mind.
Using a tolerance of $\tol = 10^{-5}$,
the integration requires \num{5015} RHS evaluations and includes many rejected steps.
Applying instead the optimized coefficients $\beta = (0.28, -0.23)$ derived
later in this manuscript results in only \num{4119} RHS evaluations and
a nearly identical final error.  A significant
performance gain is obtained by applying appropriate controller parameters, \cf
Table~\ref{tab:importance-of-controller}.

\begin{table}[!htb]
\centering\small
\caption{Performance of different controllers for \RK[BS]{5}{7}[][FSAL]:
         Number of function evaluations (\#FE) and rejected steps (\#R)
         for the linear advection problem with uniform grid using polynomials
         of degree $p = 4$ as in Section~\ref{sec:cfl-vs-error}.}
\label{tab:importance-of-controller}
\setlength{\tabcolsep}{0.75ex}
\begin{tabular*}{\linewidth}{@{\extracolsep{\fill}}c *2c r@{\hskip 0.5ex}rr@{\hskip 0.5ex}r@{\hskip 1ex}cr@{\hskip 0.5ex}r@{\hskip 1ex}c@{}}
  \toprule
  Scheme & $\beta$ & $\tol$ & \multicolumn{1}{c}{\#FE} & \multicolumn{1}{c}{\#R} & Error \\
  \midrule

  \RK[BS]{5}{7}[][FSAL]          & $(0.70, -0.40, 0.00)$
     &  $10^{-5}$ & $   5015$ & ($  132$) & \num{9.85e-07} \\
                                 & $(0.28, -0.23, 0.00)$
     &  $10^{-5}$ & $   4119$ & ($    0$) & \num{9.79e-07} \\

  \bottomrule
\end{tabular*}
\end{table}

\begin{figure}[htb]
\centering
  \hspace*{\fill}
  \begin{subfigure}{0.5\textwidth}
  \centering
    \includegraphics[width=\textwidth]{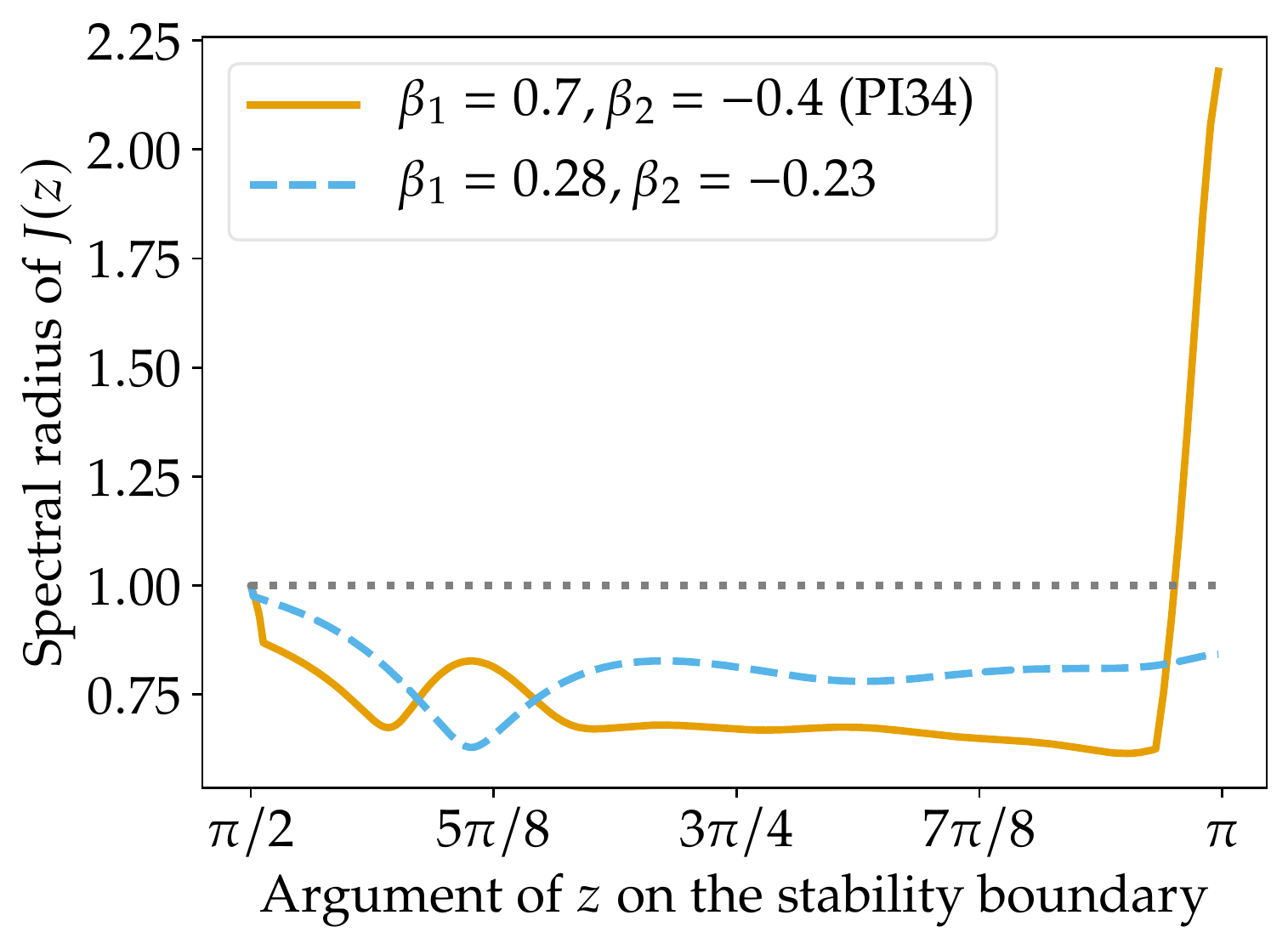}
    \caption{Spectral radius.}
  \end{subfigure}%
  \hspace*{\fill}
  \begin{subfigure}{0.31\textwidth}
  \centering
    \includegraphics[width=\textwidth]{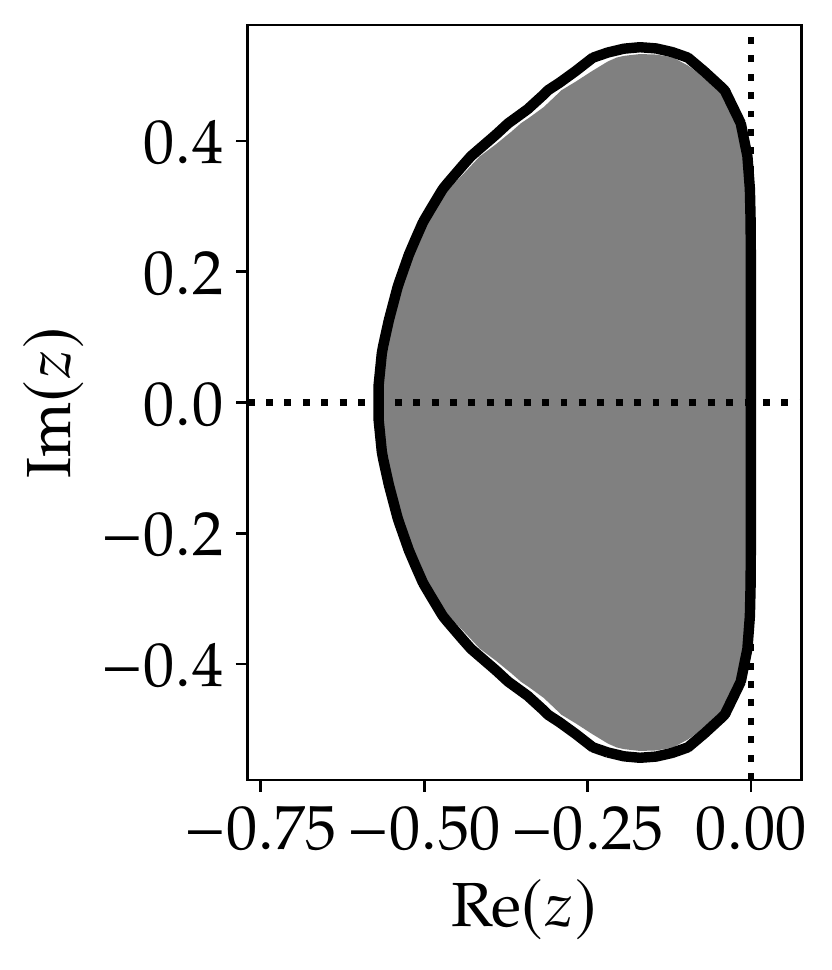}
    \caption{Stability region.}
  \end{subfigure}
  \hspace*{\fill}
  \caption{Stability region scaled by the effective number of stages and
           spectral radius of the Jacobian \eqref{eq:jacobian-control-stability}
           determining step size control stability for \RK[BS]{5}{7}[][FSAL]. The standard
           PI34 controller is unstable near the negative real axis while the
           optimized PI controller is stable along the boundary of the stability
           region.}
  \label{fig:stepsize_control_stability_bs5}
\end{figure}

The spectral radius of the Jacobian \eqref{eq:jacobian-control-stability}
determining step size control stability is plotted in
Figure~\ref{fig:stepsize_control_stability_bs5}. We see that the standard PI34
controller is unstable near the negative real axis while the optimized one
is stable.

\section{Comparison of existing methods}
\label{sec:comparing-existing-schemes}

Here, we compare some general purpose methods and schemes designed for
semidiscretizations of hyperbolic conservation laws. Since we are interested
in error-based step size control, we consider only schemes
with embedded error estimators. Hence, we consider the general purpose schemes
\begin{itemize}
  \item \RK[BS]{3}{3}[][FSAL], third-order, four-stage FSAL method of \cite{bogacki1989a32},
  \item \RK[BS]{5}{7}[][FSAL], fifth-order, eight-stage FSAL method of \cite{bogacki1996efficient},
  \item \RK[DP]{5}{6}[][FSAL], fifth-order, seven-stage FSAL method of \cite{prince1981high},
\end{itemize}
the SSP schemes
\begin{itemize}
  \item \ssp33, third-order, three-stage SSP method of
        \cite{shu1988efficient}
        with embedded method of \cite{conde2018embedded},
  \item \ssp43, third-order, four-stage SSP method of
        \cite{kraaijevanger1991contractivity}
        with embedded method of \cite{conde2018embedded} which can be implemented
        efficiently in low-storage form as described in
        Appendix~\ref{sec:ssp43},
\end{itemize}
and the low-storage methods optimized for hyperbolic conservation laws
\begin{itemize}
  \item \RK[KCL]{3}{4}[\ZRp][][C], third-order, four-stage method of \cite{kennedy2000low},
  \item \RK[KCL]{4}{5}[\ZRp][][C], fourth-order, five-stage method of \cite{kennedy2000low},
  \item \RK[KCL]{4}{5}[\ERp][][C], fourth-order, five-stage method of \cite{kennedy2000low},
  \item \RK[KCL]{5}{9}[\ZRp][][S], fifth-order, nine-stage method of \cite{kennedy2000low}.
\end{itemize}
The results shown here, obtained with three commonly-used general-purpose methods,
are typical of what we have found in tests with a much wider range of methods.
These results are sufficient to illustrate our main conclusions.
We do not consider \ssp{10}{4} of \cite{ketcheson2008highly}
with embedded method of \cite{conde2018embedded} because step size control
stability cannot be achieved for this method and any PID controller tested.
The embedded method for \ssp93 proposed in \cite{conde2018embedded} also does
not lead to step size control stability. We have created a new embedded
method with a stable optimized controller. However, it does not perform better
than \ssp43, even with manually tuned CFL-based step size control.

In the following, we will use three representative test problems to compare
the performance of these schemes. All test problems are semidiscretizations
of the compressible Euler equations in $d$ space dimensions
\begin{equation}
  \partial_t u + \sum_{j=1}^d \partial_j f^j(u) = \sigma, \quad u(0) = u_0,
\end{equation}
where the conserved variables $u = (\rho, \rho v^T, \rho e)^T$ are
the density $\rho$, the momentum $\rho v$, and the energy $\rho e$.
The flux for the spatial coordinate $j$ is
\begin{equation}
  f^j(u)
  =
  \begin{pmatrix}
    \rho v_j \\
    (\rho v v_j + p \delta_{i,j})_i \\
    (\rho e + p) v_j
  \end{pmatrix},
\end{equation}
where $p = \rho T = (\gamma-1) (\rho e - \rho v^2 / 2)$ is the pressure, $T$
the temperature, and an ideal gas law with ratio of specific heats
$\gamma = \nicefrac{7}{5}$ is assumed.
The spatial semidiscretizations use entropy-dissipative nodal DG methods with
polynomials of degree $p$ on Legendre-Gauss-Lobatto nodes with upwind
interface fluxes implemented in SSDC. We present detailed results for $p = 2$, which
is a relevant choice in practical CFD applications. The results are similar
for higher-order semidiscretizations such as polynomials of degree
$p \in \{3, 4\}$, presented in the supplementary material in more detail.

\subsection{Inviscid Taylor Green vortex}

The inviscid Taylor-Green vortex in $d=3$ space dimensions is a classical
test case to study the stability of numerical methods \cite{gassner2016split}.
The initial condition given by
\begin{equation}
\label{eq:TGV}
\begin{gathered}
  p(x) = 1 + \frac{1}{16} (\cos(2 x_1) + \cos(2 x_2)) (\cos(2 x_3) + 2),
  \\
  v_1(x) = \sin(x_1) \cos(x_2) \cos(x_3),
  \quad
  v_2(x) = -\cos(x_1) \sin(x_2) \cos(x_3),
  \\
  v_3(x) = 0,
  \quad
  T(x) = 1,
\end{gathered}
\end{equation}
with Mach number $\mach = 0.1$ is evolved in the periodic domain $[-\pi, \pi]^3$.
Unless stated otherwise, we use 8 elements per coordinate direction and a final
time of $t = 20$.
This test case is chosen as an example where the time step is mostly restricted
by stability, the solution becomes turbulent, and a relatively
low Mach number is used.

\subsection{Isentropic vortex}

The isentropic vortex is a widely used benchmark problem \cite{shu1997essentially}
with analytical solution. For the stationary case, the exact solution is given by
\begin{equation}
\label{eq:isentropic-vortex}
\begin{gathered}
  T = T_{\infty} - \frac{(\gamma-1)Ma_\infty^2\beta^2}{8 \gamma \pi^2}\exp\bigl(1-r^2\bigr),
  \\
  v_{\mathrm{t}} = \frac{r\beta}{2\pi}\exp{\left(\frac{1-r^2}{2}\right)},
  \qquad
  \rho = T^{1 / (\gamma-1)},
\end{gathered}
\end{equation}
where $r$ is the distance from the axis of the vortex and $v_{\mathrm{t}}$ is
the tangential velocity. The moving vortex solution is obtained by a uniform
translation in the direction of the velocity vector field.

Herein, the simulation domain is a cube $[-5, 5]^3$ with periodic boundaries
where the vortex rotates around the axis $(1, 1, 0)^T$, a direction not aligned
with the grid. The parameters for this test are $\gamma = 1.4$, $\mach = 0.5$,
$\beta = 5$ and $T_\infty = 1$.
Unless stated otherwise, we use 8 elements per coordinate direction for optimizing
controllers and 20 elements for examples with a final time of $t = 20$.
This test case is chosen as an example where the time step can be restricted
by accuracy for tight tolerances and because of the existence of an analytical
solution.

\subsection{Smooth flow with source terms}

The analytical solution
\begin{equation}
\label{eq:source-term}
  \rho(t,x) = \frac{3}{2} + \sin(\pi (x-t) ),
  \quad
  v(t,x) = 1,
  \quad
  p(t,x) = 1 + A_p (1 + \sin(\omega_p t)),
\end{equation}
is imposed as initial condition in the periodic domain $[-1, 1]$ and the source
term
\begin{equation}
  \sigma_{\rho e}(t,x) = \frac{A_p \omega_p}{\gamma - 1} \cos(\omega_p t)
\end{equation}
is added to the right hand side of the energy equation. The variation of the
pressure with amplitude $A_p = 50$ and frequency $\omega_p = \nicefrac{\pi}{5}$
results in a cyclic variation of the CFL restriction on the time step.
Unless stated otherwise, we use 20 elements and a final time of $t = 20$.
This test case is chosen to assess the ability of the schemes to adapt to
varying time step restrictions and because of the existence of an analytical
solution.

\subsection{Optimization of step size controllers}
\label{sec:old-controller}

As explained in Section~\ref{sec:importance-of-controller}, the choice of
appropriate step size controller parameters is important to obtain good
performance when the schemes are run at the stability limit. Hence, we
have optimized controller parameters for each scheme.
In general, the optimal time step controller parameters for a given Runge-Kutta
pair will depend somewhat on the problem under consideration. No single
controller is optimal for all test cases, but for the experiments conducted in
this work, good controllers are usually within \ca \SI{5}{\percent} of the
optimal performance.

For the low-storage schemes of \cite{kennedy2000low}, we used the
PI34 controller proposed originally with them.  We also tested the PID
controller using $\beta = (0.49, -0.34, 0.10)$ proposed in
\cite{kennedy2003additive}. We also performed an optimization of controller
parameters for each method, as follows.

We ran simulations of all three test cases described above
and measured the performance of each scheme (in terms of the number of
right-hand side evaluations).  We used a brute-force search over the domain
$\tol \in [10^{-8}, 10^{-1}]$, sampling at each power of ten, and 
$\beta_1 \in [0.1, 1.0]$, $\beta_2 \in [-0.4, -0.05]$,
$\beta_3 \in [0.0, 0.1]$ sampling at an interval of $0.01$ in each parameter,
and restricting a priori to parameter values
yielding step size control stability for the given scheme (computed using
NodePy \cite{nodepy}). The final times for these simulations were set to
$t = 8$ for \eqref{eq:TGV}, $t = 4$ for \eqref{eq:isentropic-vortex}, and
$t = 20$ for \eqref{eq:source-term} to make the brute-force optimization
feasible.
From the resulting data, consisting of thousands of runs with each
method, an overall best choice of parameters was selected as in
Section~\ref{sec:new-controller}. Usually, this kind of min-max problem was
approached by comparing the controllers minimizing the maximum, the median,
or the \SI{95}{\percent} percentile of the RHS evaluations across all CFD
simulations. Then, the final choice was made by human interaction taking
into account step size control stability and design criteria for PID controllers.

\subsection{Results for existing schemes}
\label{sec:results-existing-schemes}

For \RK[BS]{3}{3}[][FSAL], all of the controllers from Table~\ref{tab:classical-controllers}
perform reasonably well, PI42 being slightly better than the others. In general,
a wide range of controller parameters is acceptable for this scheme.
As typified by the example in Section~\ref{sec:importance-of-controller}, standard controllers
do not perform well for \RK[BS]{5}{7}[][FSAL].  We found instead that
$\beta = (0.28, -0.23, 0.00)$ is a reasonable choice for this scheme.
For \RK[DP]{5}{6}[][FSAL], the PI34 controller (which was originally designed for it by Gustafsson \cite{gustafsson1991control})
performs reasonably well in our test cases and optimized controllers like 
$\beta = (0.61, -0.27, 0.01)$ do not perform significantly better.

Subsequently, we used the optimized controller parameters and ran full simulations
(up to $t=20$) for each method with a range of tolerances.
Results are shown in Tables~\ref{tab:NFE-std-p2},
\ref{tab:NFE-ssp-p2}, and \ref{tab:NFE-ls-p2}, where polynomials of degree
$p = 2$ have been used. There, we only show results for a tolerance
$\tol = 10^{-5}$, since this choice is usually good for these small-scale test
problems. Extended details are available in the supplementary material.
For the inviscid Taylor-Green vortex, the time step is indeed restricted by stability for
most tolerances, indicated by the approximately constant number of function
evaluations, except for the very tight tolerance $\tol = 10^{-8}$ and some
schemes.
For the isentropic vortex \eqref{eq:isentropic-vortex}, the step size is
restricted by stability for tolerances $\gtrapprox 10^{-6}$.  For smaller
tolerances, the number of function evaluations increases. However, this does not result in
a significant change of the total error, which is determined mostly by the
spatial semidiscretization.
Finally, for the smooth flow with source term \eqref{eq:source-term}, the step
size is again restricted mostly by stability constraints.

\subsubsection{General purpose methods}

For very loose tolerances $\gtrapprox 10^{-3}$, \RK[BS]{5}{7}[][FSAL] and \RK[DP]{5}{6}[][FSAL] result in a
significant overhead caused by step rejections for some test cases. Otherwise,
\RK[BS]{5}{7}[][FSAL] performs better than \RK[DP]{5}{6}[][FSAL]. Other fifth-order general purpose schemes
like \RK[T]{5}{6}[][FSAL] of \cite{tsitouras2011runge} perform slightly better, usually
yielding an improvement of \ca \SI{5}{\percent}. However, \RK[BS]{3}{3}[][FSAL] is ca.
\SI{50}{\percent} more efficient as long as the time step is restricted by
stability.

These results do not change significantly if slightly higher-order
semidiscretizations are employed in space (see supplementary material),
up to polynomials of degree $p = 4$, resulting in fifth-order convergence in space.
Hence, matching the order of accuracy in space and time is not strictly necessary
if one is interested in fixed mesh sizes, especially in common CFD applications.
This remains true even if the polynomial degree is increased to $p = 7$ for
the test problems considered here. Then, the temporal error becomes significant
and the error of the fully discrete method plateaus only at relatively tight
tolerances such as $10^{-8}$. Nevertheless, \RK[BS]{3}{3}[][FSAL] is still the
most efficient method for such high-order methods and tight tolerances.

\begin{table}[htb]
\centering\small
\caption{Performance of general purpose schemes:
         Number of function evaluations (\#FE), rejected steps (\#R), and
         $L^2$ error of the density for
         the inviscid Taylor-Green vortex \eqref{eq:TGV},
         the isentropic vortex \eqref{eq:isentropic-vortex} ,
         and the flow with source term \eqref{eq:source-term}
         using polynomials of degree $p = 2$.}
\label{tab:NFE-std-p2}
\setlength{\tabcolsep}{0.75ex}
\begin{tabular*}{\linewidth}{@{\extracolsep{\fill}}c *2c r@{\hskip 0.5ex}rr@{\hskip 0.5ex}r@{\hskip 1ex}cr@{\hskip 0.5ex}r@{\hskip 1ex}c@{}}
  \toprule
   & & & \multicolumn{2}{c}{TGV \eqref{eq:TGV}} & \multicolumn{3}{c}{Isent. vortex \eqref{eq:isentropic-vortex}} & \multicolumn{3}{c}{Source term \eqref{eq:source-term}} \\
  Scheme & $\beta$ & $\tol$ & \multicolumn{1}{c}{\#FE} & \multicolumn{1}{c}{\#R} & \multicolumn{1}{c}{\#FE} & \multicolumn{1}{c}{\#R} & Error & \multicolumn{1}{c}{\#FE} & \multicolumn{1}{c}{\#R} & Error \\
  \midrule

  \RK[BS]{3}{3}[][FSAL]          & $(0.60, -0.20, 0.00)$
     &  $10^{-5}$ & $   5256$ & ($    1$) & $   1476$ & ($    0$) & \num{5.77e-04} & $  20682$ & ($    4$) & \num{1.77e-03} \\
  \RK[BS]{5}{7}[][FSAL]          & $(0.28, -0.23, 0.00)$
     &  $10^{-5}$ & $   7722$ & ($    1$) & $   2159$ & ($    0$) & \num{5.78e-04} & $  30356$ & ($    2$) & \num{1.77e-03} \\
  \RK[DP]{5}{6}[][FSAL]          & $(0.70, -0.40, 0.00)$
     &  $10^{-5}$ & $   8014$ & ($    3$) & $   2217$ & ($    0$) & \num{5.78e-04} & $  31430$ & ($    4$) & \num{1.77e-03} \\

  \bottomrule
\end{tabular*}
\end{table}

\subsubsection{SSP methods}

The popular method \ssp33 can be equipped with the PI34 controller to give
acceptable step size control performance; slightly better behavior can be
achieved by choosing $\beta = (0.70, -0.37, 0.05)$. For loose and medium
tolerances, this scheme performs similarly to \RK[BS]{3}{3}[][FSAL].
\RK[BS]{3}{3}[][FSAL] is significantly more efficient than \ssp33 at tight tolerances.

\ssp43 performs \ca \SI{50}{\percent} better than \ssp33 or \RK[BS]{3}{3}[][FSAL] at loose
and medium tolerances for the inviscid Taylor-Green vortex using the
optimized controller $\beta = (0.55, -0.27, 0.05)$. At loose tolerances,
it is also \ca \SI{15}{\percent} more efficient than \RK[BS]{3}{3}[][FSAL] for the
isentropic vortex. However, the number of RHS evaluations increases
drastically for tighter tolerances, making \ssp43 less efficient than
\RK[BS]{3}{3}[][FSAL] for these parameters.
The results for the flow with source term are
similar but less pronounced. Hence, \ssp43 can be more efficient than the
best schemes so far but the embedded method does not seem to be reliable
enough to make the choice of the tolerance as robust as for other schemes.
Additionally, the choice of appropriate controller parameters can be crucial
for \ssp43, since some standard controllers do not perform well.

For higher polynomial degrees $p \in \{3, 4\}$, \ssp43 is still a very interesting
method that can even beat \RK[BS]{3}{3}[][FSAL]. \ssp33 is less efficient than
\ssp43 also for these higher polynomial degrees. For even higher polynomial
degrees such as $p = 7$, the situation changes a bit since the temporal error
becomes significant. While \ssp43 is still the most efficient method so far
for medium tolerances, it becomes less efficient than \RK[BS]{3}{3}[][FSAL]
for the vortex problems at a tolerance of $10^{-8}$, since it is less optimized
for accuracy than that general purpose method.

\begin{table}[htb]
\centering\small
\caption{Performance of SSP schemes:
         Number of function evaluations (\#FE), rejected steps (\#R), and
         $L^2$ error of the density for
         the inviscid Taylor-Green vortex \eqref{eq:TGV},
         the isentropic vortex \eqref{eq:isentropic-vortex},
         and the flow with source term \eqref{eq:source-term}
         using polynomials of degree $p = 2$.}
\label{tab:NFE-ssp-p2}
\setlength{\tabcolsep}{0.75ex}
\begin{tabular*}{\linewidth}{@{\extracolsep{\fill}}c *2c r@{\hskip 0.5ex}rr@{\hskip 0.5ex}r@{\hskip 1ex}cr@{\hskip 0.5ex}r@{\hskip 1ex}c@{}}
  \toprule
   & & & \multicolumn{2}{c}{TGV \eqref{eq:TGV}} & \multicolumn{3}{c}{Isent. vortex \eqref{eq:isentropic-vortex}} & \multicolumn{3}{c}{Source term \eqref{eq:source-term}} \\
  Scheme & $\beta$ & $\tol$ & \multicolumn{1}{c}{\#FE} & \multicolumn{1}{c}{\#R} & \multicolumn{1}{c}{\#FE} & \multicolumn{1}{c}{\#R} & Error & \multicolumn{1}{c}{\#FE} & \multicolumn{1}{c}{\#R} & Error \\
  \midrule

  \ssp33                         & $(0.70, -0.37, 0.05)$
     &  $10^{-5}$ & $   5261$ & ($    3$) & $   1646$ & ($    0$) & \num{5.77e-04} & $  20681$ & ($    3$) & \num{1.77e-03} \\
  \ssp43                         & $(0.55, -0.27, 0.05)$
     &  $10^{-5}$ & $   3438$ & ($    1$) & $   1902$ & ($    0$) & \num{5.78e-04} & $  20638$ & ($    3$) & \num{1.77e-03} \\

  \bottomrule
\end{tabular*}
\end{table}

\subsubsection{Low-storage methods}

Some standard controllers like PI34 perform mostly acceptably well for \RK[KCL]{3}{4}[\ZRp][][C]
(based on the number of step rejections). Nevertheless, an optimized controller
with parameters $\beta = (0.50, -0.35, 0.10)$ results in a few percent fewer
function evaluations.
However, \RK[BS]{3}{3}[][FSAL] is up to \SI{20}{\percent} more efficient, in accordance with
the real stability interval scaled by the effective number of stages, which is
three for \RK[BS]{3}{3}[][FSAL] because of the FSAL property.

The PI34 controller does not perform well for the other low-storage schemes.
For \RK[KCL]{4}{5}[\ERp][][C] (but not for the other schemes), the PID controller with
$\beta = (0.49, -0.34, 0.10)$ proposed in \cite{kennedy2003additive} performs
much better.
An optimized controller with parameters $\beta = (0.41, -0.28, 0.08)$ performs
even slightly better, making this scheme more efficient than \RK[BS]{3}{3}[][FSAL] for the
inviscid Taylor-Green vortex and slightly more efficient for the isentropic vortex.
However, \RK[BS]{3}{3}[][FSAL] is still better
for the other test case.
\RK[KCL]{4}{5}[\ZRp][][C] and \RK[KCL]{5}{9}[\ZRp][][S] were more challenging
for optimizing controller parameters and less efficient than \RK[KCL]{4}{5}[\ERp][][C].

As for $p = 2$, \RK[KCL]{4}{5}[\ERp][][C] is usually the most efficient existing
low-storage method of \cite{kennedy2000low} for $p \in \{3, 4\}$, which can
also be more efficient than \RK[BS]{3}{3}[][FSAL] for the vortex problems.
However, \ssp43 is even more efficient there. Additionally, the sensitivity of
the step size controller for the low-storage methods is bigger than for
\RK[BS]{3}{3}[][FSAL].
For $p = 7$, all of the low-storage methods considered here result in a
non-negligible amount of step rejections for the inviscid Taylor-Green vortex.
Nevertheless, the fourth-order accurate methods can be up to \SI{15}{\percent}
more efficient than \RK[BS]{3}{3}[][FSAL] there. Nevertheless, \ssp43 is still
more efficient for this test problem. At tight tolerances such as $10^{-8}$,
\RK[KCL]{4}{5}[\ERp][][C] is the most efficient method considered so far for
the isentropic vortex and the smooth flow with source term.

\begin{table}[htb]
\centering\small
\caption{Performance of low-storage schemes:
         Number of function evaluations (\#FE), rejected steps (\#R), and
         $L^2$ error of the density for
         the inviscid Taylor-Green vortex \eqref{eq:TGV},
         the isentropic vortex \eqref{eq:isentropic-vortex},
         and the flow with source term \eqref{eq:source-term}
         using polynomials of degree $p = 2$.}
\label{tab:NFE-ls-p2}
\setlength{\tabcolsep}{0.75ex}
\begin{tabular*}{\linewidth}{@{\extracolsep{\fill}}c *2c r@{\hskip 0.5ex}rr@{\hskip 0.5ex}r@{\hskip 1ex}cr@{\hskip 0.5ex}r@{\hskip 1ex}c@{}}
  \toprule
   & & & \multicolumn{2}{c}{TGV \eqref{eq:TGV}} & \multicolumn{3}{c}{Isent. vortex \eqref{eq:isentropic-vortex}} & \multicolumn{3}{c}{Source term \eqref{eq:source-term}} \\
  Scheme & $\beta$ & $\tol$ & \multicolumn{1}{c}{\#FE} & \multicolumn{1}{c}{\#R} & \multicolumn{1}{c}{\#FE} & \multicolumn{1}{c}{\#R} & Error & \multicolumn{1}{c}{\#FE} & \multicolumn{1}{c}{\#R} & Error \\
  \midrule

  \RK[KCL]{3}{4}[\ZRp][][C]      & $(0.50, -0.35, 0.10)$
     &  $10^{-5}$ & $   6326$ & ($    1$) & $   1874$ & ($    0$) & \num{5.78e-04} & $  24886$ & ($    3$) & \num{1.77e-03} \\
  \RK[KCL]{4}{5}[\ZRp][][C]      & $(0.29, -0.24, 0.02)$
     &  $10^{-5}$ & $   4897$ & ($   65$) & $   1522$ & ($    0$) & \num{5.79e-04} & $  23902$ & ($    5$) & \num{1.77e-03} \\
  \RK[KCL]{4}{5}[\ERp][][C]      & $(0.41, -0.28, 0.08)$
     &  $10^{-5}$ & $   4587$ & ($   42$) & $   1472$ & ($    0$) & \num{5.79e-04} & $  23817$ & ($    4$) & \num{1.77e-03} \\
  \RK[KCL]{5}{9}[\ZRp][][S]      & $(0.49, -0.34, 0.10)$
     &  $10^{-5}$ & $   6320$ & ($    3$) & $   2108$ & ($    0$) & \num{5.78e-04} & $  30374$ & ($    3$) & \num{1.77e-03} \\

  \bottomrule
\end{tabular*}
\end{table}

\subsubsection{Discussion}

All of the general purpose schemes make use of the FSAL technique. Additionally,
the stability region of the embedded scheme is always at least as big as the
one of the main method.
Although the 2R low-storage schemes were optimized for convection-dominated
problems, they were outperformed for all test problems and at almost all tolerances
by the general purpose method
\RK[BS]{3}{3}[][FSAL]. Possible reasons for this are that the 2R method
coefficients are chosen subject to more stringent low-storage requirements, they
do not exploit the FSAL technique, and they have embedded methods with a
stability region that is smaller
than that of the main method in some areas. When the time step is restricted by
stability, this can effectively reduce the allowable time step for the main method.

\begin{figure}[!htb]
\centering
  \begin{subfigure}{0.33\textwidth}
    \includegraphics[width=\textwidth]{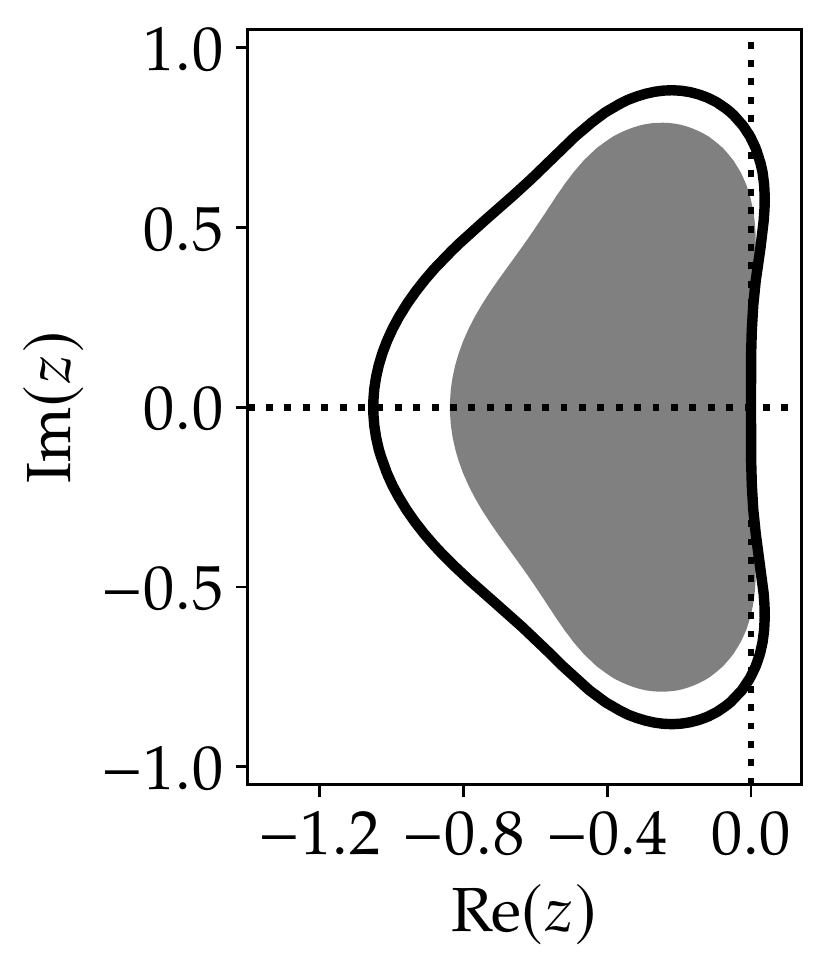}
    \caption{\RK[BS]{3}{3}[][FSAL].}
  \end{subfigure}%
  \hspace*{\fill}
  \begin{subfigure}{0.33\textwidth}
    \includegraphics[width=\textwidth]{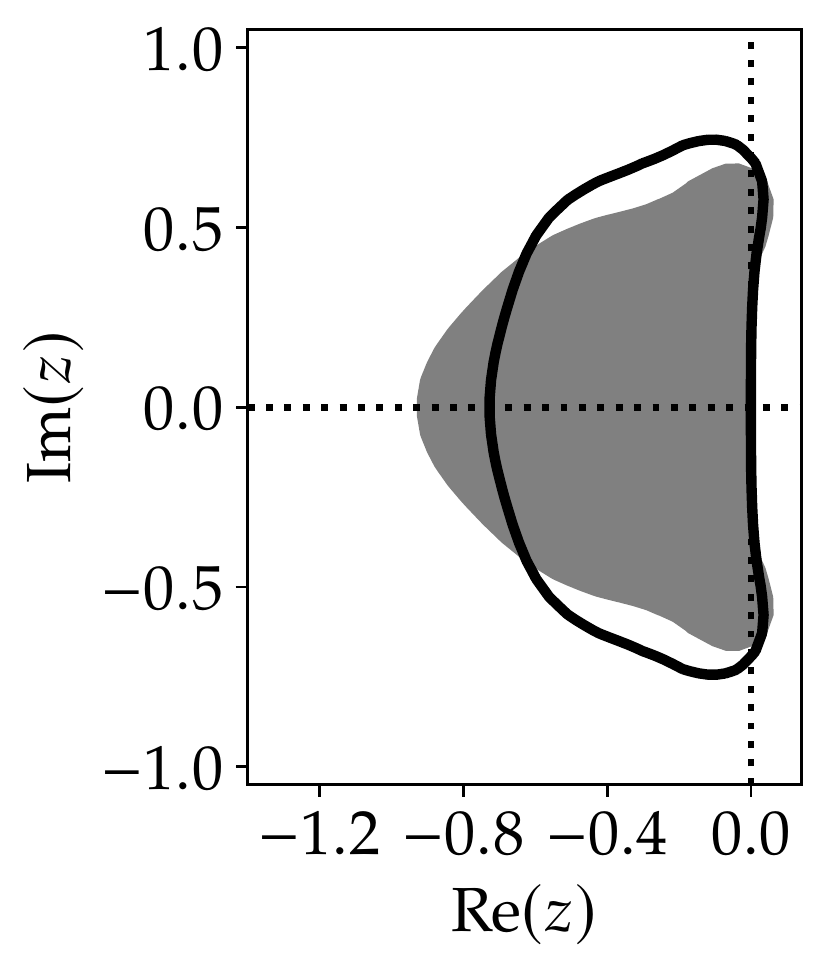}
    \caption{\RK[KCL]{4}{5}[\ERp][][C].}
  \end{subfigure}%
  \hspace*{\fill}
  \begin{subfigure}{0.33\textwidth}
    \includegraphics[width=\textwidth]{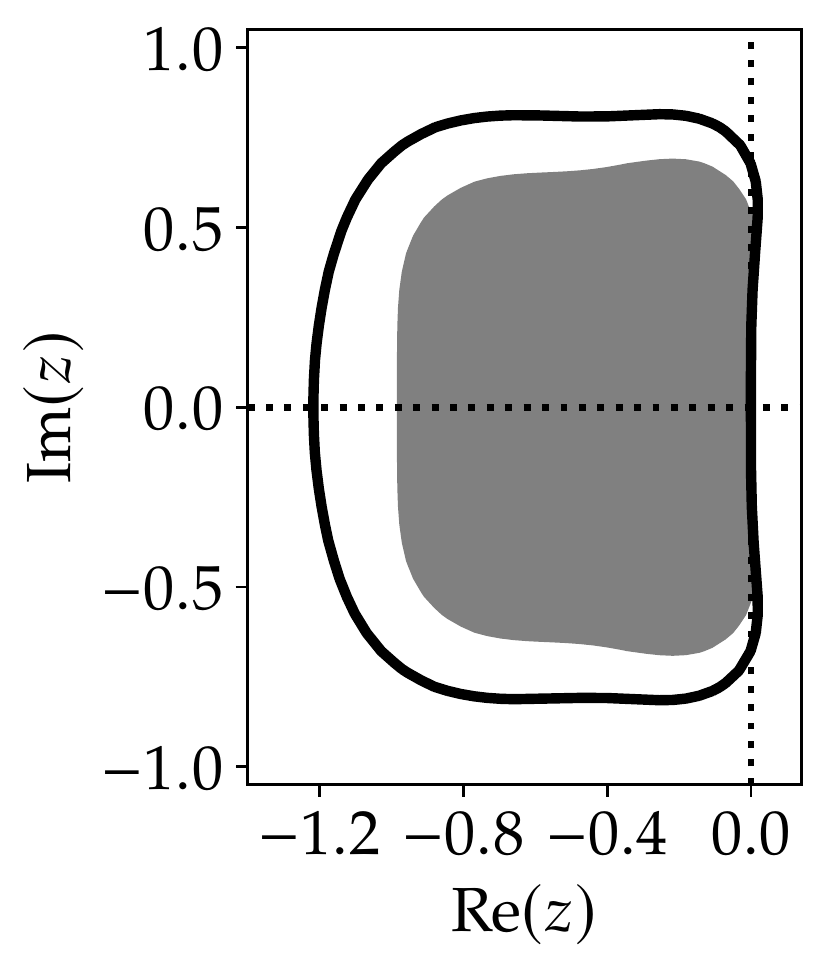}
    \caption{\RK{3}{5}[\ESstarp][FSAL].\label{fig:stability_regionsc}}
  \end{subfigure}%
  \caption{Stability regions scaled by the effective number of stages of three
           representative methods (taking the FSAL property into account).
           The stability region of the main method is marked in gray and the
           boundary of the embedded method's stability region is drawn as
           black line.}
  \label{fig:stability_regions}
\end{figure}

This last point is illustrated in Figure~\ref{fig:stability_regions}, which
shows the stability regions (scaled by the effective number of stages)
of the main and embedded method for three pairs.
We see that although the stability region of \RK[BS]{3}{3}[][FSAL] includes
less of the real axis than that of \RK[KCL]{4}{5}[\ERp][][C], the embedded
method for \RK[BS]{3}{3}[][FSAL] extends further than that of \RK[KCL]{4}{5}[\ERp][][C].
The last stability region in the figure corresponds to a new method developed in
the next section.  Like \RK[BS]{3}{3}[][FSAL], it has the useful property that
the stability region of the embedded method contains that of the main method.

In the results described above, the behavior of the methods for the inviscid Taylor-Green
vortex was often slightly different than for the other test cases. This can
partly be explained by the lower Mach number chosen for this example. Indeed,
numerical experiments show that low Mach numbers put more stress on real axis
stability than on the rest of the spectrum generated by linear advection.
Hence, methods with stability regions that include more parts of the negative
real axis but are not optimal for the linear advection spectrum can perform
better for low Mach numbers; see also \cite{citro2020optimal}.
In this article, we focus on applications in CFD with medium to high Mach
numbers.
However, flows with small Mach numbers are usually computed using incompressible
solvers and implicit time integration methods. Hence, we do not focus
on this regime in this article. Nevertheless, we use the inviscid Taylor-Green
vortex as test case to study the step size control stability on the negative real
axis.

\section{New Optimized Runge-Kutta pairs}
\label{sec:optimized-RK}

In the previous section we developed optimized step size controllers for
existing Runge-Kutta pairs.  Now we consider the optimization of
Runge-Kutta pairs themselves (along with controllers).
To do so, we begin with the \ESstarp methods of
\cite{parsani2013optimized}, without embedded error estimator.
Then, we design an embedded method and optimized controller parameters.
The embedded method is optimized for step size control stability, good error
metrics (see \cite{kennedy2000low} and the supplementary material for the present
work), to have a large stability region that includes that of the main method,
and to have coefficients that are not too large.
The resulting schemes given by double precision floating point numbers are
optimized further using extended precision numbers in Julia \cite{bezanson2017julia}
and the package Optim.jl \cite{mogensen2018optim}, such that the order conditions
are satisfied at least to quadruple precision.
Coefficients of the new optimized methods are available in the accompanying
repository \cite{ranocha2021optimizedCoefficients} in full precision. Double
precision coefficients are given in Appendix~\ref{sec:coefficients}.
The stability region of a representative method is shown in
Figure~\ref{fig:stability_regionsc}.

We also developed new pairs from scratch, based on the approach
used in \cite{parsani2013optimized}.
Specifically, we compute the Fourier footprint of the spectral
element semidiscretization of the linear advection equation by varying the
direction of the wave propagation velocity vector, the solution orientation, and
the wave vector module and construct
optimized stability polynomials using the algorithm described in \cite{ketcheson2013optimal}.
Afterwards, low-storage Runge-Kutta schemes are constructed by minimizing their
principal error constants, given their class, the number of stages $s$, the
order of accuracy $q$, and the optimized stability polynomial as constraints.
These optimizations are carried out using RK-Opt \cite{RKopt}, based on the
optimization toolbox of MATLAB.  However, the resulting methods did not perform
better than the pairs based on starting with methods from \cite{parsani2013optimized}.

\subsection{Optimization of controller parameters}
\label{sec:new-controller}

The controller parameters for these new pairs are optimized using the same
approach as described in Section~\ref{sec:comparing-existing-schemes}.

\begin{figure}[!htb]
\centering
  \begin{subfigure}{0.49\textwidth}
    \includegraphics[width=\textwidth]{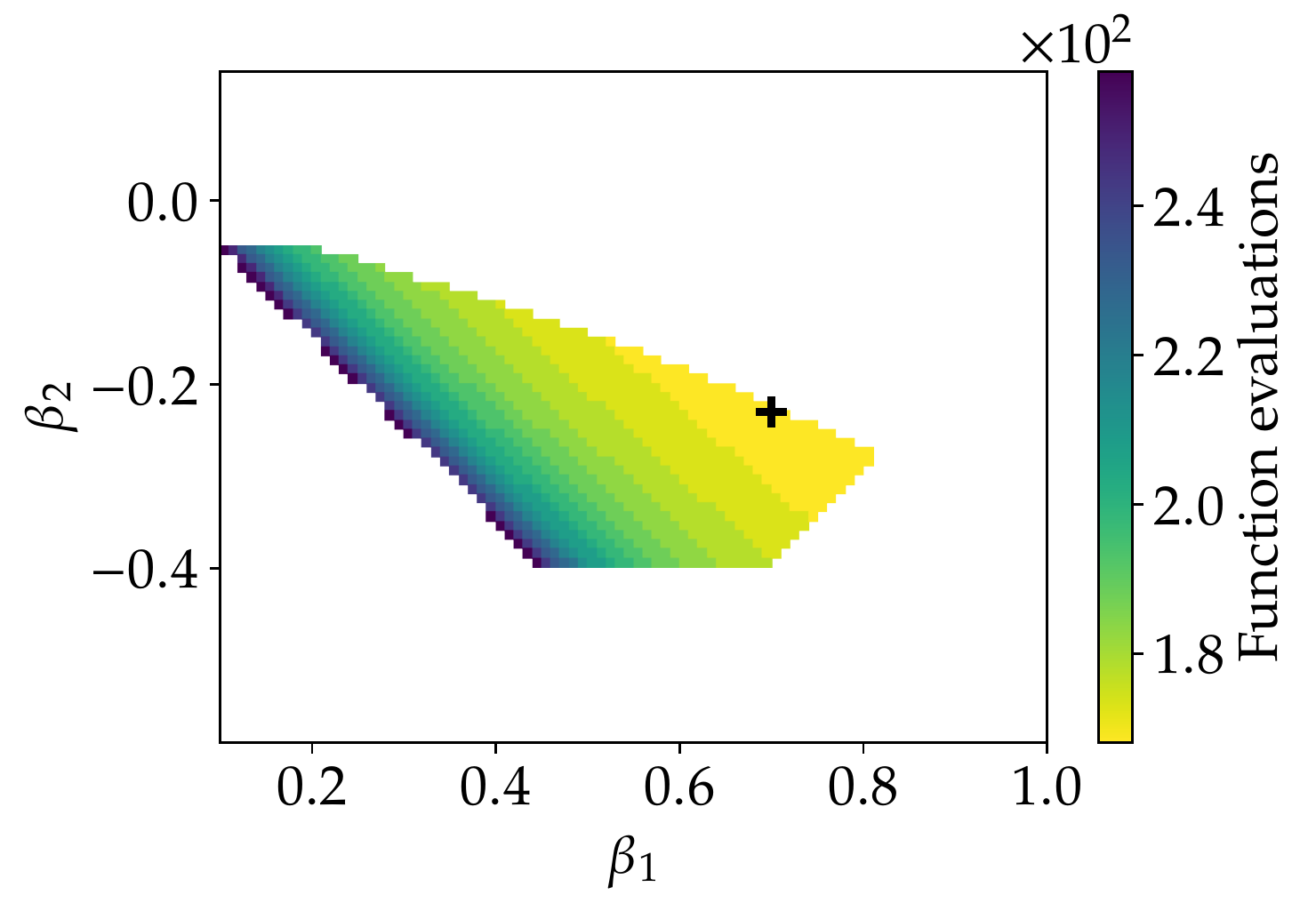}
    \caption{Isentropic vortex \eqref{eq:isentropic-vortex},
             the chosen controller is within \SI{3.1}{\percent} of the minimal
             \#FE.}
  \end{subfigure}%
  \hspace*{\fill}
  \begin{subfigure}{0.49\textwidth}
    \includegraphics[width=\textwidth]{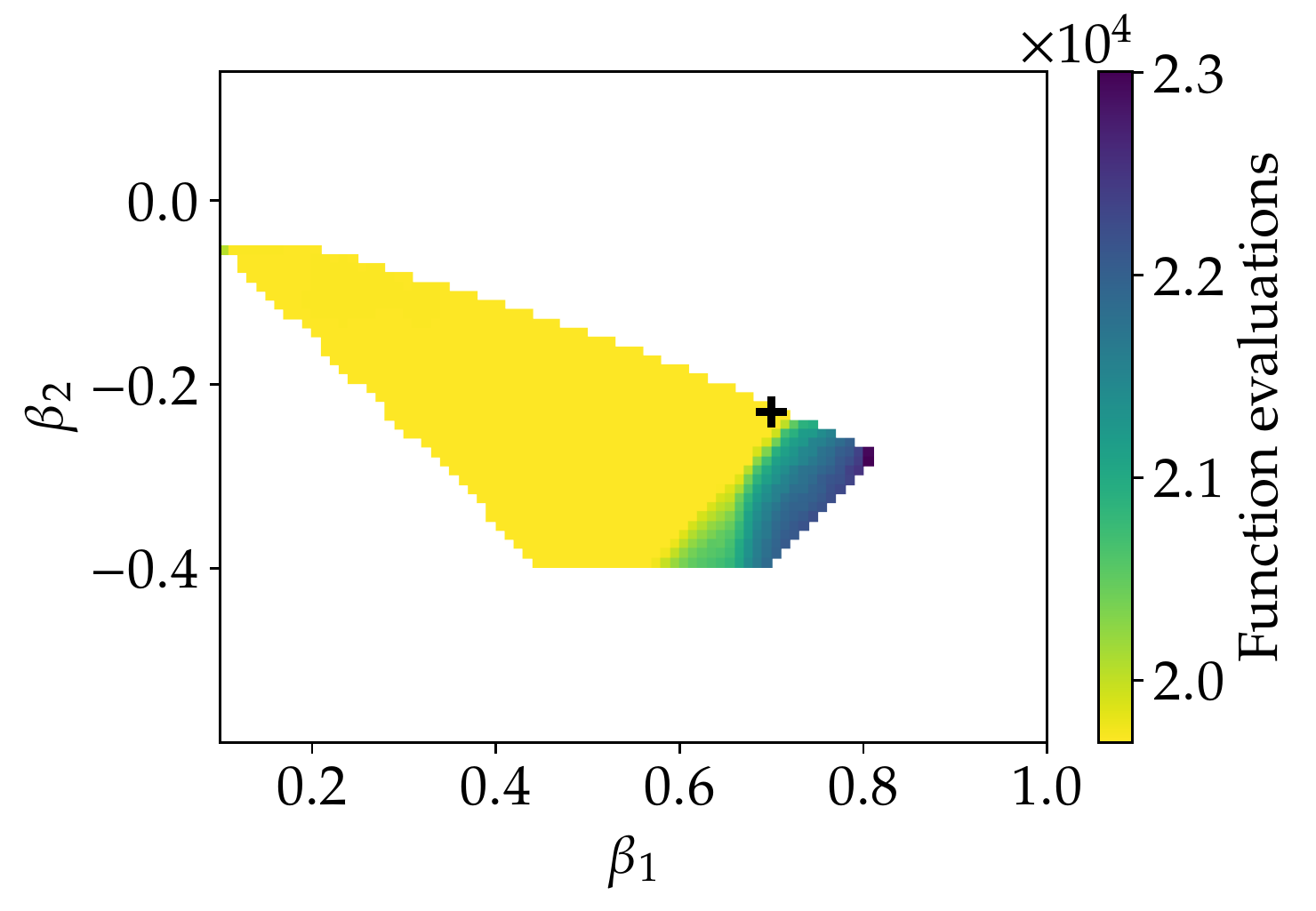}
    \caption{Flow with source term \eqref{eq:source-term},
             the chosen controller is within \SI{0.1}{\percent} of the minimal
             \#FE.}
  \end{subfigure}%
  \caption{Performance of different controllers with $\beta_3 = 0$ for
           \RK{3}{5}[\ESstarp][FSAL] and two of the test problems with tolerance
           $\tol = 10^{-5}$. The chosen controller for this scheme uses
           $\beta = (0.70, -0.23, 0.00)$ (marked with a black +).
           The number of function evaluations
           (\#FE) is visualized only for those controllers that result in step size
           control stability along the whole boundary of the main method's
           stability region.}
  \label{fig:controller_performance_3SstarpFSAL}
\end{figure}

Typical performance results of the optimization procedure are shown in
Figure~\ref{fig:controller_performance_3SstarpFSAL} using \RK{3}{5}[\ESstarp][FSAL]
as an example.
For the isentropic vortex \eqref{eq:isentropic-vortex} with $\tol = 10^{-5}$,
the temporal accuracy starts to play a role and the controllers are not
necessarily limited by stability. In this regime, controllers with larger
$\beta_1$ and $\beta_2$ closer to zero perform better; they are more near
to the simple deadbeat (I-) controller, which is in some sense optimal in
the asymptotic regime.
In contrast, the controllers operate near the stability boundary for the
test case with source term \eqref{eq:source-term}. Here, controllers with more
damping and smaller $\beta_1$ perform better.
To find an acceptable controller for the scheme \RK{3}{5}[\ESstarp][FSAL], both kinds
of problems have to be considered, seeking a compromise between efficiency in
the asymptotic regime and near the stability boundary.

\subsection{Results for new schemes}

Analogously to Tables~\ref{tab:NFE-std-p2}--\ref{tab:NFE-ls-p2}, results
are summarized in Table~\ref{tab:NFE-opt-p2} for the new optimized low-storage
schemes with error control for $p = 2$; extended details and results for
higher-order spatial semidiscretizations using solution polynomials of degree
$p \in \{3, 4, 7\}$ are available in the supplementary material.

\begin{table}[htb]
\centering\small
\caption{Performance of optimized low-storage schemes:
         Number of function evaluations (\#FE), rejected steps (\#R), and
         $L^2$ error of the density for
         the inviscid Taylor-Green vortex \eqref{eq:TGV},
         the isentropic vortex \eqref{eq:isentropic-vortex},
         and the flow with source term \eqref{eq:source-term}
         using polynomials of degree $p = 2$.}
\label{tab:NFE-opt-p2}
\setlength{\tabcolsep}{0.75ex}
\begin{tabular*}{\linewidth}{@{\extracolsep{\fill}}c *2c r@{\hskip 0.5ex}rr@{\hskip 0.5ex}r@{\hskip 1ex}cr@{\hskip 0.5ex}r@{\hskip 1ex}c@{}}
  \toprule
   & & & \multicolumn{2}{c}{TGV \eqref{eq:TGV}} & \multicolumn{3}{c}{Isent. vortex \eqref{eq:isentropic-vortex}} & \multicolumn{3}{c}{Source term \eqref{eq:source-term}} \\
  Scheme & $\beta$ & $\tol$ & \multicolumn{1}{c}{\#FE} & \multicolumn{1}{c}{\#R} & \multicolumn{1}{c}{\#FE} & \multicolumn{1}{c}{\#R} & Error & \multicolumn{1}{c}{\#FE} & \multicolumn{1}{c}{\#R} & Error \\
  \midrule

  \RK{3}{5}[\ESstarp]            & $(0.64, -0.31, 0.04)$
     &  $10^{-5}$ & $   4467$ & ($    2$) & $   1262$ & ($    0$) & \num{5.77e-04} & $  19692$ & ($    3$) & \num{1.77e-03} \\
  \RK{3}{5}[\ESstarp][FSAL]      & $(0.70, -0.23, 0.00)$
     &  $10^{-5}$ & $   4469$ & ($    1$) & $   1308$ & ($    0$) & \num{5.77e-04} & $  19690$ & ($    2$) & \num{1.77e-03} \\
  \RK{4}{9}[\ESstarp]            & $(0.25, -0.12, 0.00)$
     &  $10^{-5}$ & $   4205$ & ($    1$) & $   1478$ & ($    0$) & \num{5.79e-04} & $  18992$ & ($    1$) & \num{1.78e-03} \\
  \RK{4}{9}[\ESstarp][FSAL]      & $(0.38, -0.18, 0.01)$
     &  $10^{-5}$ & $   4207$ & ($    1$) & $   1371$ & ($    0$) & \num{5.80e-04} & $  18984$ & ($    1$) & \num{1.77e-03} \\
  \RK{5}{10}[\ESstarp]           & $(0.47, -0.20, 0.06)$
     &  $10^{-5}$ & $   5372$ & ($    2$) & $   1512$ & ($    0$) & \num{5.78e-04} & $  24107$ & ($    3$) & \num{1.78e-03} \\
  \RK{5}{10}[\ESstarp][FSAL]     & $(0.45, -0.13, 0.00)$
     &  $10^{-5}$ & $   5364$ & ($    1$) & $   1533$ & ($    0$) & \num{5.78e-04} & $  24112$ & ($    2$) & \num{1.77e-03} \\

  \bottomrule
\end{tabular*}
\end{table}

In general, the novel schemes are more efficient than all methods tested
in Section \ref{sec:comparing-existing-schemes}. In particular, the novel
third-order schemes are up to \SI{18}{\percent}
more efficient than \RK[BS]{3}{3}[][FSAL], in accordance with the relative
lengths of the real stability intervals.
They are also up to \SI{5}{\percent} more efficient than \RK[KCL]{4}{5}[\ERp][][C] with the
optimized PID controller for the inviscid Taylor-Green vortex and up to
\SI{13}{\percent} more efficient for the isentropic vortex. Recall that
\RK[BS]{3}{3}[][FSAL] is more efficient than \RK[KCL]{4}{5}[\ERp][][C] for the
other test problem.
Only \ssp43 is more efficient than \RK{3}{5}[\ESstarp][FSAL] for the inviscid Taylor-Green
vortex, in accordance with the particularly large real stability interval,
\cf Section~\ref{sec:results-existing-schemes}. However, the new schemes are more efficient
at realistic (medium to high) Mach numbers, for which they have been optimized.

The optimized fourth-order schemes can be even more efficient for the inviscid Taylor-Green
vortex and the smooth flow with source term at medium tolerances, giving an
improvement of a few percent. For the isentropic vortex, the third-order schemes
are still up to \SI{6}{\percent} more efficient.
The optimized fourth-order schemes are up to \SI{25}{\percent} more efficient
than the best corresponding methods of \cite{kennedy2000low} with the controller
recommended there.
However, the fourth-order accurate main methods of \cite{parsani2013optimized}
make it particularly difficult to design good embedded methods and controllers.
This can already be seen in the optimized controller coefficients, where
the magnitudes of $\beta_1$ and $\beta_2$ differ less than for other optimized
methods. While the controllers can be tuned to result in acceptable performance
for these test problems, they do not necessarily lead to good performance for
other setups.

The optimized fifth-order schemes are less efficient than the optimized third-
and fourth-order schemes for these test cases unless the tolerance is very
tight (so the spatial error dominates and the influence of the time
integrator is negligible).
These schemes (used with the controllers designed here) are much more efficient
than the corresponding method of \cite{kennedy2000low} (used with the controller
prescribed there), by up to \SI{25}{\percent} for the source term problem, up to
\SI{35}{\percent} for the isentropic vortex, and up to \SI{18}{\percent} for the
inviscid Taylor-Green vortex.

For $p = 3$, the third-order accurate schemes are the most efficient ones of
the optimized low-storage methods, except for the inviscid Taylor-Green vortex,
where the fourth-order methods are up to \SI{6}{\percent} more efficient.
For $p = 4$, the fourth-order accurate schemes are the most efficient new ones
for these experiments, followed closely by the third-order accurate methods.
For $p = 7$, the fourth-order accurate schemes are still the most efficient new
ones; the third-order methods do not match the same small errors for tight
tolerances and their temporal error dominates the spatial one. However, the
fourth-order methods are difficult to control for loose tolerances, resulting in
a significant number of step rejections. For sufficiently tight tolerances,
the optimized fourth-order methods are more efficient than the fourth-order
methods of \cite{kennedy2000low}.

In general, FSAL methods are often more efficient than non-FSAL schemes,
especially at loose tolerances. Thus, we recommend to use the novel
\ESstarp\ (FSAL) methods for hyperbolic problems where the time step is restricted
mostly by stability constraints.

Optimization of Runge-Kutta pairs for higher numbers of stages was not as
successful. While we were able to obtain schemes with good theoretical
properties, their performance did not show improvement compared to the schemes
listed above. Improvements of the underlying optimization algorithms or
the imposition of additional constraints might lead to better schemes
in the future. However, the novel schemes developed in this
work are already a significant improvement over the state of the art and
perform well.

\subsection{Further optimizations}

As shown in \cite{parsani2013optimized} and the numerical experiments above,
the spatial error usually dominates the temporal error. Hence, it is interesting
to optimize lower-order time integrators for higher-order spatial discretizations.
Such an approach is presented in \cite{parsani2012optimized}, but focused on
first-order accurate time integrators. The results shown there demonstrate some
speedup compared to the schemes presented in \cite{parsani2013optimized}, but
these come with a reduced accuracy.

\begin{figure}[htb]
\centering
  \includegraphics[width=0.3\textwidth]{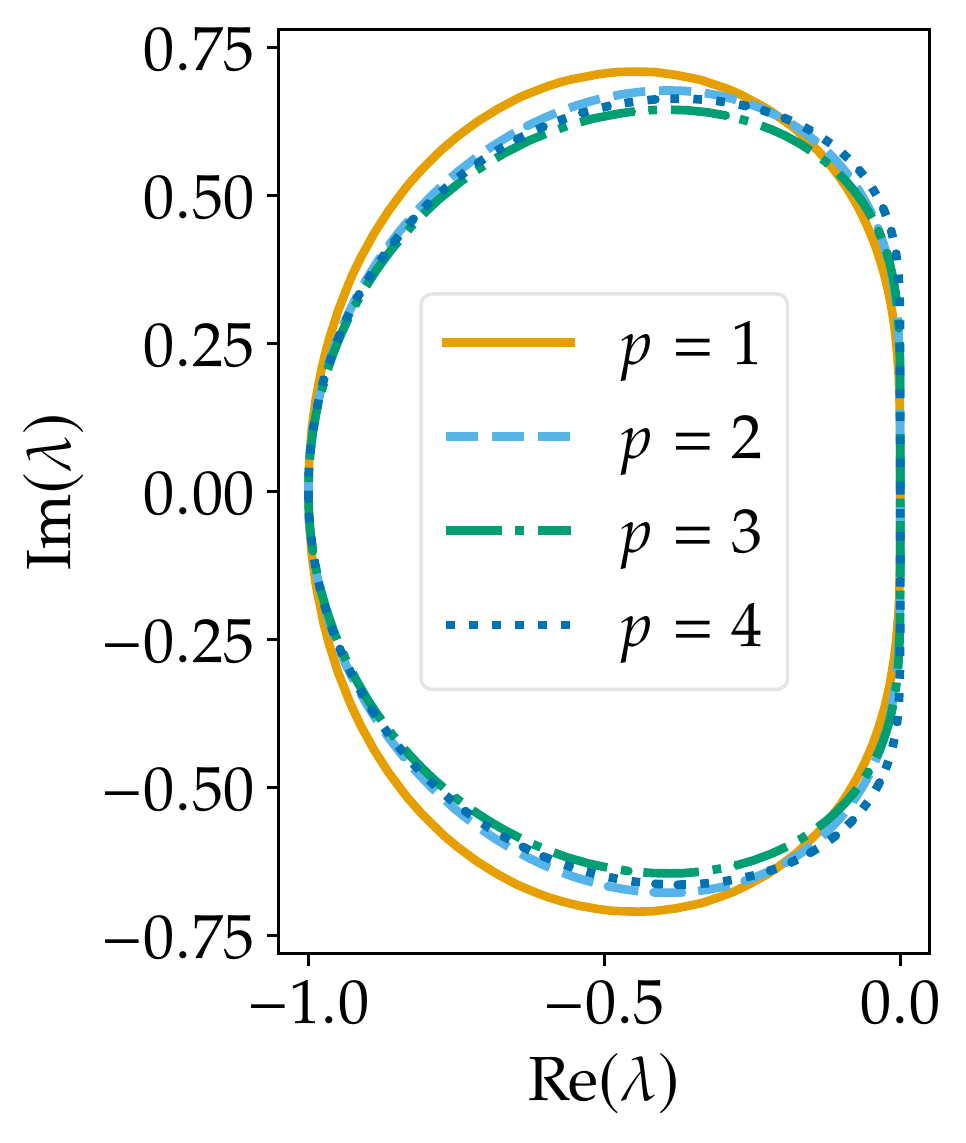}
  \caption{Convex hulls of the spectra of the spectral element semidiscretizations
           used for the optimization of stability polynomials as in
           \cite{parsani2013optimized}. The spectra are scaled such that
           $\min \mathrm{Re}(\lambda) = -1$.}
  \label{fig:spectra_SD}
\end{figure}

Here, we choose third-order accurate Runge-Kutta methods and optimize them for
a fifth-order spatial discretization. This resulted in a speedup compared to
the optimized schemes described in the previous sections for some test cases;
for other test cases, no speedup could be observed. This is in accordance with
the general similarity of the scaled convex hulls of the spectra for different
polynomial degrees $p \in {1, 2, 3, 4}$ shown in Figure~\ref{fig:spectra_SD}.
Hence, we do not pursue this path of research further.

\section{Additional numerical experiments and comparisons}
\label{sec:numerical-experiments}

Hitherto, a careful selection of test cases was used to demonstrate issues and
design criteria for explicit Runge-Kutta schemes applied to semidiscretizations
of hyperbolic conservation laws. Next, more involved examples are used to
demonstrate that the novel methods can be applied successfully to large-scale
CFD problems including the compressible Navier-Stokes equations.

To be useful for engineering and applied problems in CFD, a CFL-based control
must be automated as much as possible. Therefore, we use the approach described
in Section~\ref{sec:cfl-vs-error} also for the viscous CFL number.
The normalizing factor $\sigma$ in \eqref{eq:cfl-factor} is chosen depending on
the solution polynomial degree such that a method with a real stability interval of $2$
is stable for the linear advection-diffusion equation on a uniform grid with a
CFL factor $\cfl = 1$. On top of that, a safety factor of $0.95$ is applied,
\cf \cite{jahdali2021optimized}.

Except for the viscous shock described next, the other simulations conducted
here start from a checkpoint of a developed solution and run on 8 nodes of Shaheen XC40 using
32 CPU cores each (one compute node of Shaheen XC40\footnote{Shaheen XC40 is the petascale supercomputer hosted at KAUST, which features 6{,}174 dual
socket compute nodes based on 16 core Intel Haswell processors running at \SI{2.3}{GHz}.
Each node is equipped with \SI{128}{GB} of DDR4 memory running at \SI{2.3}{GHz}. Overall, the system has a
total of 197{,}568 processor cores and \SI{790}{TB} of aggregate memory.}).
The general purpose and SSP methods are implemented using the explicit
Runge-Kutta interface of PETSc. The other methods are implemented using their
respective low-storage forms in PETSc.

\subsection{Viscous shock}
\label{sec:VS}

The propagating viscous shock is a classical test problem for the compressible Navier-Stokes
equations. The momentum $\fnc{V}$ of the analytical solution satisfies the ODE
\begin{equation}
  \alpha \fnc{V}\frac{\partial\fnc{V}}{\partial x}
  - (\fnc{V}-1)(\fnc{V}-\fnc{V}_{f}) = 0,
  \qquad -\infty \le x \le +\infty.
\end{equation}
The solution of this ODE can be written implicitly as
\begin{equation}
\label{eq:VS}
  x-\frac{1}{2}\alpha\left(\log\left|(\fnc{V}(x_1)-1)(\fnc{V}(x_1)-\fnc{V}_{f})\right|+\frac{1+\fnc{V}_{f}}{1-\fnc{V}_{f}}\log\left|\frac{\fnc{V}(x_1)-1}{\fnc{V}(x_1)-\fnc{V}_{f}}\right|\right) = 0,
\end{equation}
where
\begin{equation}
  \fnc{V}_{f} \equiv \frac{\fnc{U}_{L}}{\fnc{U}_{R}}, \qquad
  \alpha \equiv \frac{2\gamma}{\gamma + 1}\frac{\,\mu}{Pr\dot{\fnc{M}}}.
\end{equation}
Here, $\fnc{U}_{L/R}$ are the known velocities to the left and right of the
shock at $\pm\infty$, $\dot{\fnc{M}}$ is the constant mass flow across the shock,
$Pr = 3/4$ is the Prandtl number, and $\mu$ is the dynamic viscosity. The mass
and total enthalpy are constant across the shock. Moreover, the momentum and
energy equations become redundant.

For our tests, we compute $\fnc{V}$ from \eqref{eq:VS}
to machine precision using bisection.
The moving shock solution is obtained by applying a uniform translation to
the above solution. Initially, at $t = 0$, the shock is located at the center
of the domain. We use the parameters $\mach=2.5$, $\reynolds=10$,
and $\gamma=1.4$ in the domain given by $x \in [-0.5,0.5]$ till the final time
$t = 2$.
The boundary conditions are prescribed by penalizing the numerical solution
against the analytical solution, which is also used to prescribe the
initial condition.

\begin{table}[htb]
\centering\small
\caption{Number of function evaluations (\#FE), rejected steps (\#R), and
         $L^2$ error of the density for
         the viscous shock \eqref{eq:VS}
         using polynomials of degree $p$.}
\label{tab:NFE-VS}
\setlength{\tabcolsep}{0.75ex}
\begin{tabular*}{\linewidth}{@{\extracolsep{\fill}}c *2c r@{\hskip 0.5ex}r@{\hskip 1ex}cr@{\hskip 0.5ex}r@{\hskip 1ex}c@{}}
  \toprule
   & & & \multicolumn{3}{c}{$p = 2$} & \multicolumn{3}{c}{$p = 4$} \\
  Scheme & $\beta$ & $\tol$ & \multicolumn{1}{c}{\#FE} & \multicolumn{1}{c}{\#R} & Error & \multicolumn{1}{c}{\#FE} & \multicolumn{1}{c}{\#R} & Error \\
  \midrule

  \RK[BS]{3}{3}[][FSAL]          & $(0.60, -0.20, 0.00)$
     &  $10^{-5}$ & $    615$ & ($    1$) & \num{3.76e-03} & $   3842$ & ($    3$) & \num{5.42e-05} \\
  \ssp43                         & $(0.55, -0.27, 0.05)$
     &  $10^{-5}$ & $    450$ & ($    1$) & \num{3.76e-03} & $   2502$ & ($    3$) & \num{5.42e-05} \\
  \RK[KCL]{4}{5}[\ERp][][C]      & $(0.41, -0.28, 0.08)$
     &  $10^{-5}$ & $    552$ & ($    2$) & \num{3.76e-03} & $   3342$ & ($   29$) & \num{5.42e-05} \\
  \RK{3}{5}[\ESstarp][FSAL]      & $(0.70, -0.23, 0.00)$
     &  $10^{-5}$ & $    533$ & ($    0$) & \num{3.76e-03} & $   3270$ & ($    2$) & \num{5.42e-05} \\
  \RK{4}{9}[\ESstarp]            & $(0.25, -0.12, 0.00)$
     &  $10^{-5}$ & $    542$ & ($    0$) & \num{3.76e-03} & $   3403$ & ($    7$) & \num{5.47e-05} \\
  \RK{5}{10}[\ESstarp][FSAL]     & $(0.45, -0.13, 0.00)$
     &  $10^{-5}$ & $    643$ & ($    0$) & \num{3.76e-03} & $   3931$ & ($    3$) & \num{5.42e-05} \\

  \bottomrule
\end{tabular*}
\end{table}

Some results for the most promising methods and optimized controllers are shown in
Table~\ref{tab:NFE-VS}; extended details are available in the supplementary
material.
The new scheme \RK{3}{5}[\ESstarp][FSAL] is \ca \SI{18}{\percent} more
efficient than \RK[BS]{3}{3}[][FSAL] for relevant tolerances, in accordance with
the relative real stability intervals.
\ssp43 is a very promising scheme for this kind of problem because of its
improved stability properties around the negative real axis. In particular,
\ssp43 is \ca \SI{50}{\percent} more efficient than \RK[BS]{3}{3}[][FSAL]
for relevant tolerances, also in accordance with the relative real stability
intervals.
Except for very tight tolerances and low solution polynomial degrees, the step size
controllers detect the stability constraint accurately; the spatial error
dominates and the error of the time integration schemes is negligible.

\subsection{NASA juncture flow}

We consider the NASA juncture flow problem as described in
\cite[Section~3.8]{parsani2021ssdc}. The NASA juncture flow test was designed
to validate CFD for wing juncture trailing edge separation and progression,
and it is a collaborative effort between CFD computationalists and
experimentalists \cite{rumsey2018nasajuncture}. 
Specifically, the NASA juncture flow experiment is a series of wind tunnel 
tests conducted in the NASA Langley subsonic tunnel
to collect
validation data in the juncture region of a wing-body configuration \cite{rumsey2016nasajuncture}.

\begin{figure}[htbp!]
\centering
\includegraphics[width=0.95\columnwidth]{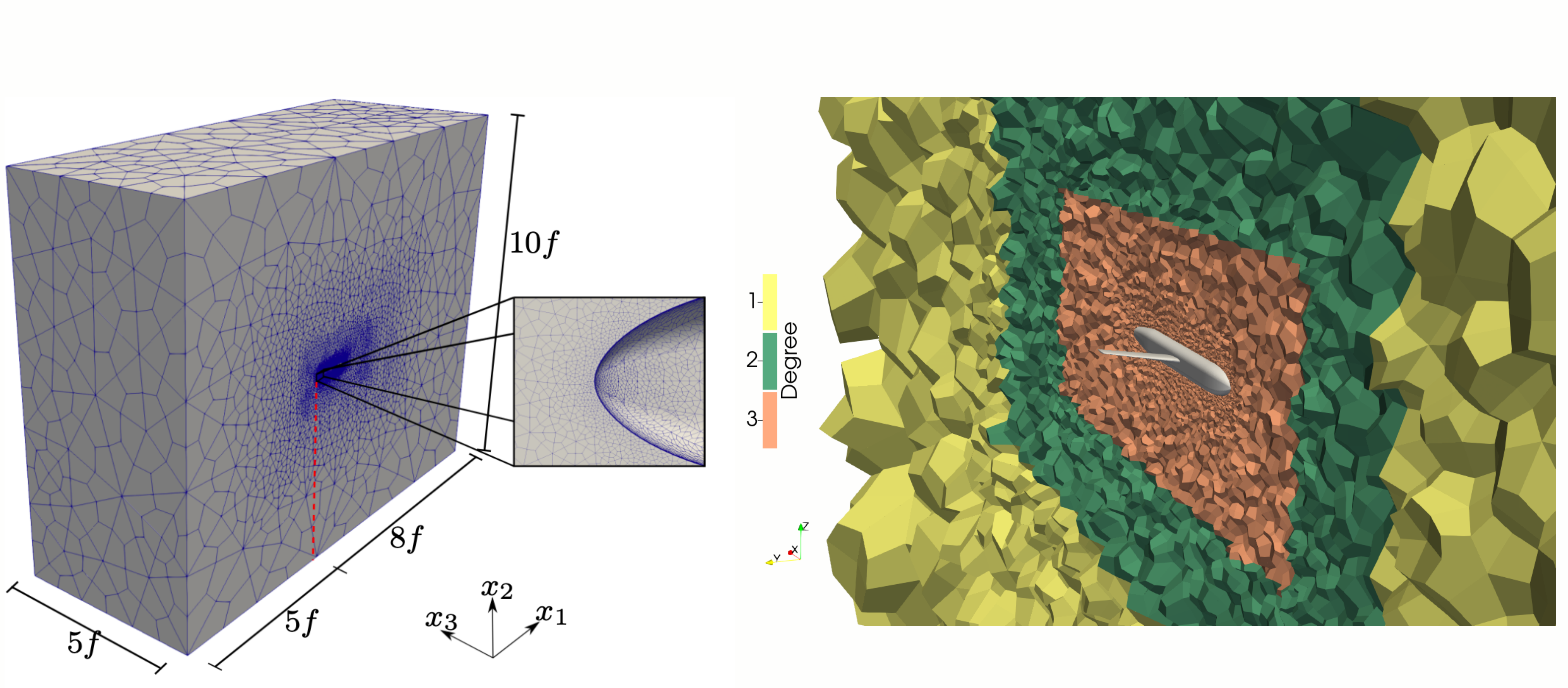}
\caption{Solution polynomial degree distribution, computational domain and boundary mesh elements for the NASA juncture experiment 
	\cite{parsani2021ssdc}; $f$ is the fuselage length.}
\label{fig:junc_all}
\end{figure}

Here, we simulate the NASA juncture flow with a wing based on the DLR-F6 geometry and a leading edge
horn to mitigate the effect of the horseshoe vortex over the wing-fuselage juncture \cite{langtry2005drag}. 
A general view of the geometry is shown in Figure~\ref{fig:junc_all}(b).
The model crank chord is $\ell=\SI{557.1}{mm}$, the
wing span is $77.89\ell$, and the fuselage length is $f=8.69\ell$. 
The wing leading edge horn meets the fuselage at $x_1=3.45\ell$, and the wing root 
trailing-edge is located at $x_1=5.31\ell$.
In the wind tunnel, the model is mounted on a sting aligned with the fuselage axis. 
The sting is attached to a mast that emerges from the wind tunnel floor.
The Reynolds number is $\reynolds=2.4 \times 10^6$ and the freestream Mach number is
$\mach=0.189$. The angle of attack is AoA~$=-2.5^\circ$. We perform simulations in
free air conditions, ignoring both the sting and the mast.

As shown in Figure~\ref{fig:junc_all}(b), the grid is subdivided into three blocks, corresponding to three
different approximation degrees, $p$, for the solution field.
In particular, we use $p=1$ in the far-field region, $p=3$ in the region surrounding the  
model, and $p=2$ elsewhere. In total, we use $\approx 6.762 \times 10^5$
hexahedral elements and $\approx 4.091 \times 10^7$ degrees of freedom (DOFs).
We highlight that the boundary layer thickness over the fuselage for $x =1,000$--$2,000$~mm is
about \SI{16}{mm}, while it is about \SI{20}{mm} over
the wing upstream of the separation bubble \cite{iyer2020nasajunction}. 
In the present simulation we use between eight and nine solution points in the boundary layer thickness
$\delta_{99}$. The mesh features a maximum aspect ratio of \ca~$110$.
The grid is constructed using the commercial software Pointwise~V18.3 released in September~2019; solid boundaries are described using a quadratic mesh.

\begin{figure}[!htb]
\centering
  \begin{subfigure}{0.52\textwidth}
    \includegraphics[width=\textwidth]{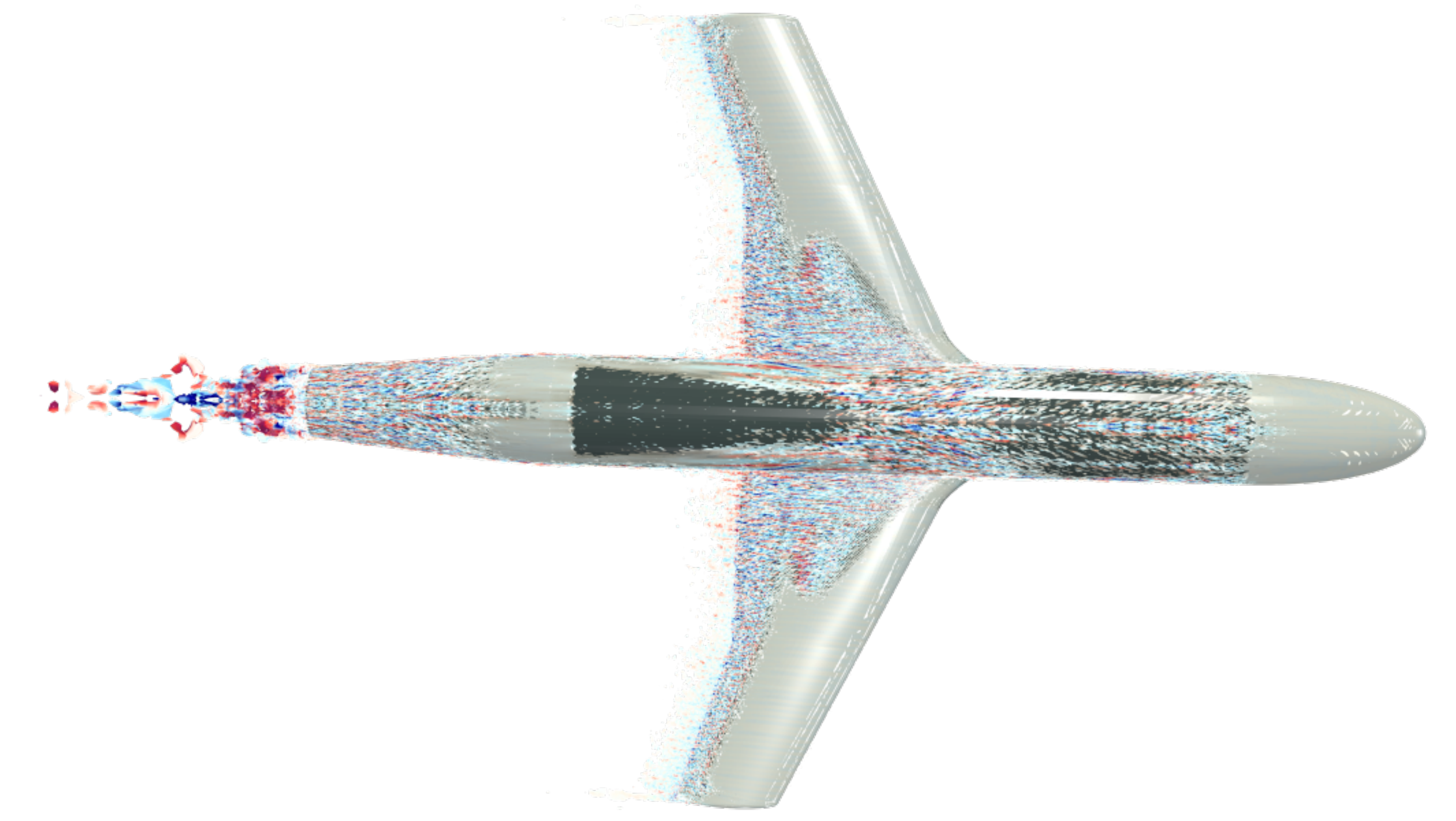}
    \caption{Top view.}
  \end{subfigure}%
  \hspace*{\fill}
  \begin{subfigure}{0.49\textwidth}
    \includegraphics[width=\textwidth]{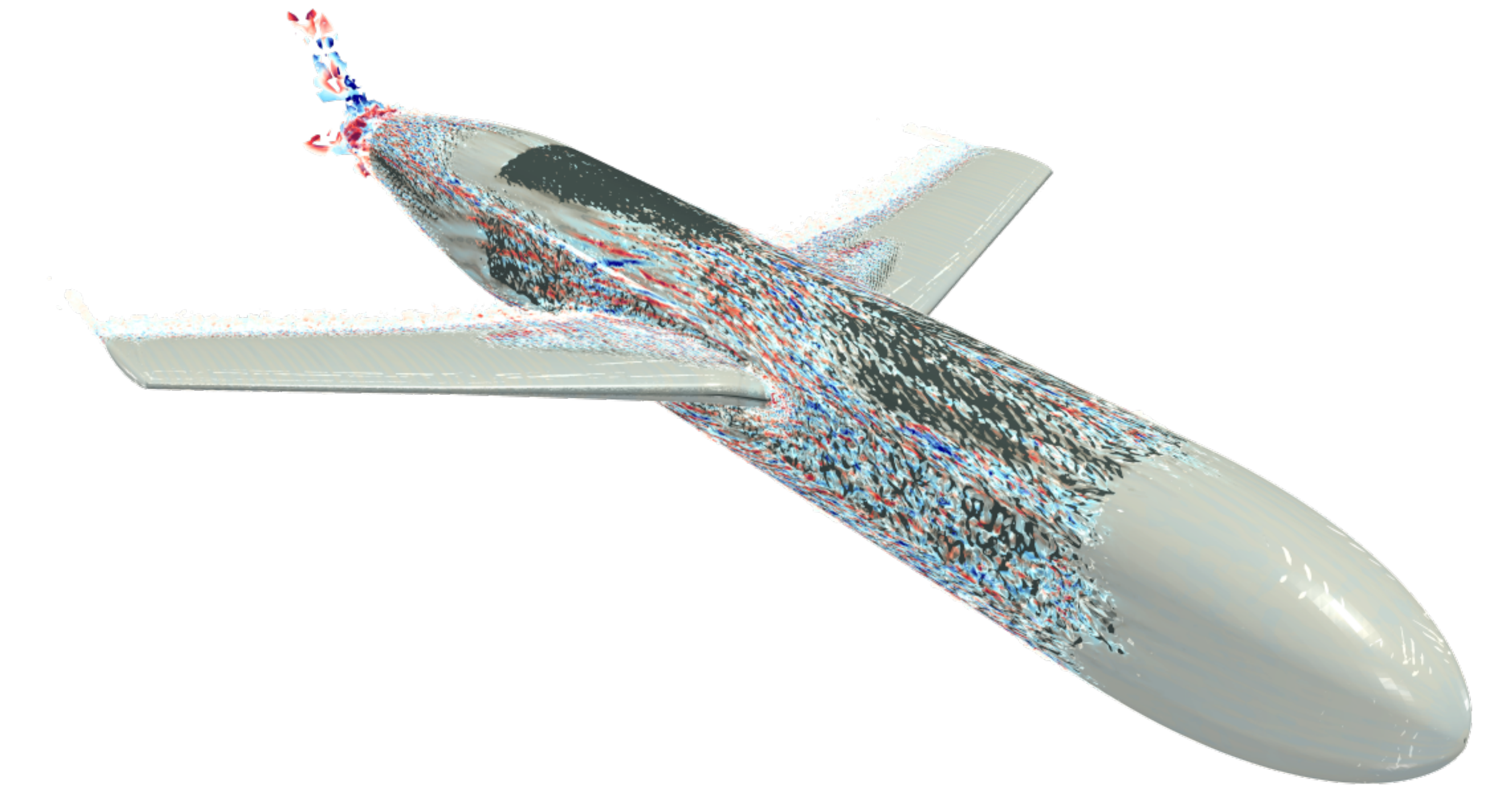}
    \caption{3D view.}
  \end{subfigure}%
  \caption{$Q$-criterion colored by the velocity magnitude of the NASA juncture flow.}
  \label{fig:qcrit_nasa_juncture}
\end{figure}

Figure~\ref{fig:qcrit_nasa_juncture} shows the $Q$-criterion colored by the velocity magnitude of flow
past the aircraft. The separation of the flow on the wing near the junction with the fuselage is visible.

\begin{table}[!htb]
\sisetup{output-exponent-marker=}
\sisetup{scientific-notation=fixed, fixed-exponent=0}
\centering\small
\caption{Performance of different methods for error and CFL-based step size
         controllers:
         Number of function evaluations (\#FE), rejected steps (\#R), and
         wall-clock time in seconds for the NASA juncture flow problem.}
\label{tab:NASA-juncture}
\setlength{\tabcolsep}{0.75ex}
\begin{tabular*}{\linewidth}{@{\extracolsep{\fill}}c *2c r@{\hskip 0.5ex}rr@{\hskip 0.5ex}r@{\hskip 1ex}cr@{\hskip 0.5ex}r@{\hskip 1ex}c@{}}
  \toprule
  Scheme & $\beta$ & $\tol$/$\cfl$ & \multicolumn{1}{c}{\#FE} & \multicolumn{1}{c}{\#R} & Wall-clock time (s) \\
  \midrule

  \RK[BS]{3}{3}[][FSAL]          & $(0.60, -0.20, 0.00)$
     &  $\tol = 10^{-8}$ & $   1786$ & $    1$ & \tablenum{5.8649e+02} \\
   & &  $\cfl =  1.0   $ & $   3598$ &         & \tablenum{1.17160e+03} \smallskip\\
  \ssp43                         & $(0.55, -0.27, 0.05)$
     &  $\tol = 10^{-8}$ & $   1172$ & $    1$ & \tablenum{3.7618e+02} \\
   & &  $\cfl =   1.0  $ & $   2340$ &         & \tablenum{7.4745e+02} \smallskip\\
  \RK{3}{5}[\ESstarp][FSAL]      & $(0.70, -0.23, 0.00)$
     &  $\tol = 10^{-8}$ & $   1527$ & $    1$ & \tablenum{4.9658e+02} \\
   & &  $\cfl =   1.0  $ & $   3056$ &         & \tablenum{9.9135e+02} \smallskip\\
  \RK{4}{9}[\ESstarp]            & $(0.25, -0.12, 0.00)$
     &  $\tol = 10^{-8}$ & $   1467$ & $    1$ & \tablenum{4.6920e+02} \\
   & &  $\cfl =   1.0  $ & $   2862$ &         & \tablenum{9.1217e+02} \smallskip\\
  \RK{5}{10}[\ESstarp][FSAL]     & $(0.45, -0.13, 0.00)$
     &  $\tol = 10^{-8}$ & $   1864$ & $    3$ & \tablenum{5.9706e+02} \\
   & &  $\cfl =   1.0  $ & $   3661$ &         & \tablenum{1.17010e+03} \\

  \bottomrule
\end{tabular*}
\end{table}

A summary of the performance of the different methods is presented in
Table~\ref{tab:NASA-juncture}. Here, the CFL adaptor with $\cfl = 1.0$
tuned for linear advection-diffusion works for all RK methods. However, it was
significantly less efficient than the error-based controller with a conservative tolerance
of $10^{-8}$; the CFL controller used \ca \SI{50}{\percent} more RHS evaluations
and wall-clock time. Thus, a tedious manual tuning to increase the CFL factor would be
necessary to match the efficiency of the error-based controller which just
works out of the box.

\RK[BS]{3}{3}[][FSAL] is an efficient general purpose method for this CFD problem.
Nevertheless, the optimized third- and fourth-order accurate methods are more
efficient. Interestingly, \ssp43 is again significantly more efficient, nearly
\SI{50}{\percent} faster than \RK[BS]{3}{3}[][FSAL].

\subsection{Viscous flow past a Formula 1 front wing}

\begin{figure}[!ht]
\centering
\includegraphics[width=0.98\columnwidth]{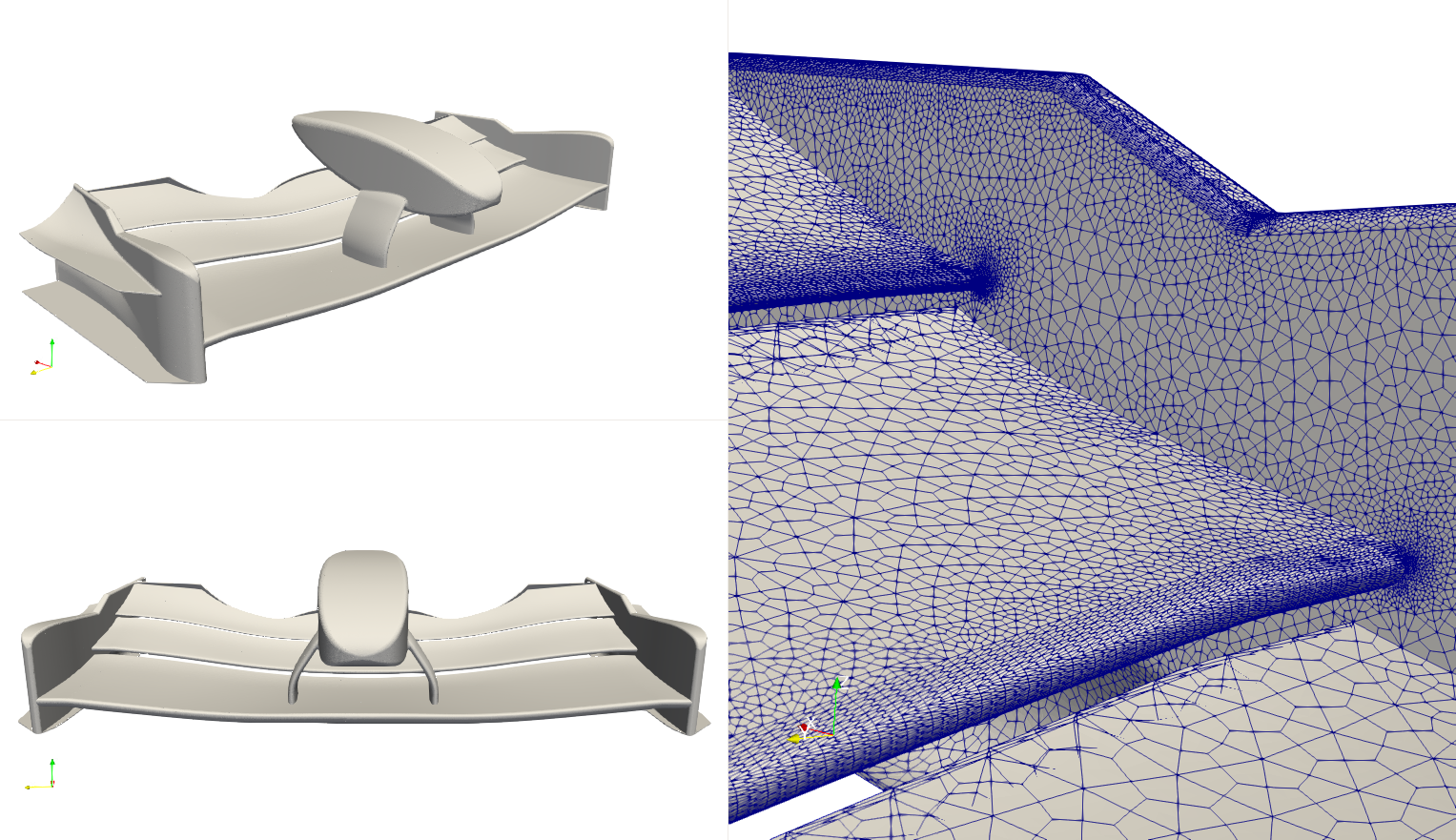}
\caption{Overview of the Imperial Front Wing.}
\label{fig:ifw}
\end{figure}

Here, we consider the flow past a Formula 1 front wing with a relatively complex geometry, supported by the
availability of a CAD model and experimental results \cite{pegrum2007experimental}. We refer to this test case as the Imperial Front Wing,
originally based on the front wing and endplate design of the McLaren 17D race car \cite{buscariolo2019}. The panel of Figure~\ref{fig:ifw}
gives an overview of the Imperial Front Wing geometry.
We denote by $h$ the distance between the ground and the lowest part of the front wing endplate and by $c$ the 
chord length of the main element. The position of the wing in the tunnel is further characterized by a
pitch angle of $1.094^{\circ}$. Here we use $h/c = 0.36$ which can be considered as a relatively low front ride height, with high
ground effect and hence higher loads on the wing. The corresponding Reynolds number is $\reynolds = 2.2 \times 10^5$,
based on the main element chord $c$ of \SI{250}{mm} and a free stream velocity $U$ of \SI{25}{m/s}. The Mach number is set to
$\mach = 0.036$. This corresponds to a practically incompressible flow.

The computational domain is divided into $3.4 \times 10^6$ hexahedral elements with a maximum aspect ratio of \ca~$250$. Two different semidiscretizations
with solution polynomials of degree $p = 1$ and $p = 2$ are used. The grid is constructed using the commercial software Pointwise~V18.3 released in September~2019; solid boundaries are described using a quadratic mesh.

\begin{figure}[!htb]
\centering
  \begin{subfigure}{0.49\textwidth}
    \includegraphics[width=\textwidth]{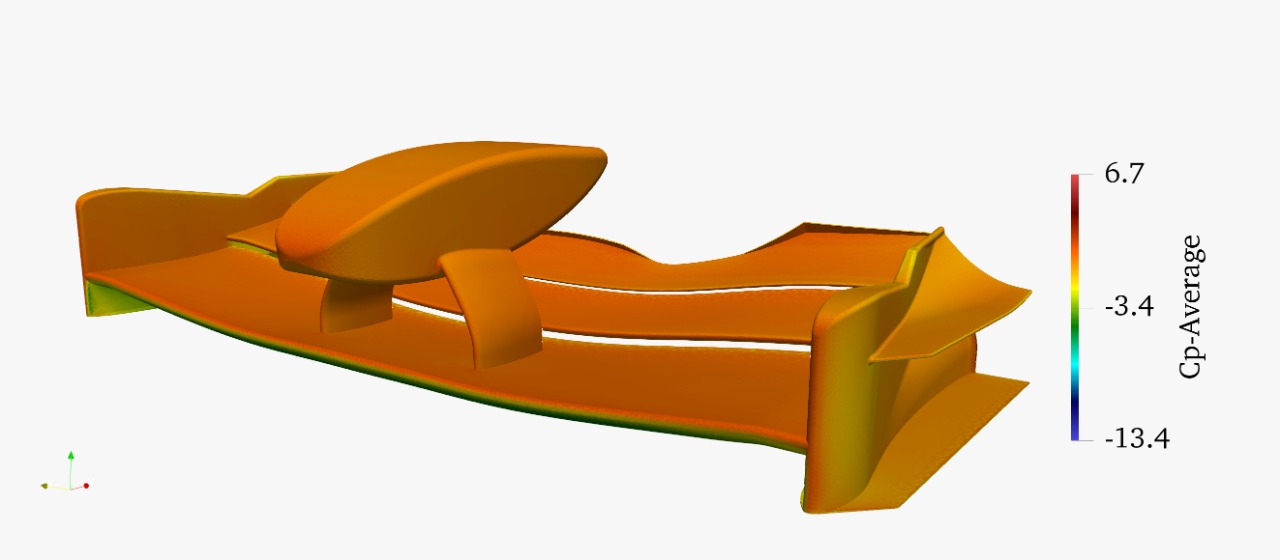}
    \caption{Top.}
  \end{subfigure}%
  \hspace*{\fill}
  \begin{subfigure}{0.49\textwidth}
    \includegraphics[width=\textwidth]{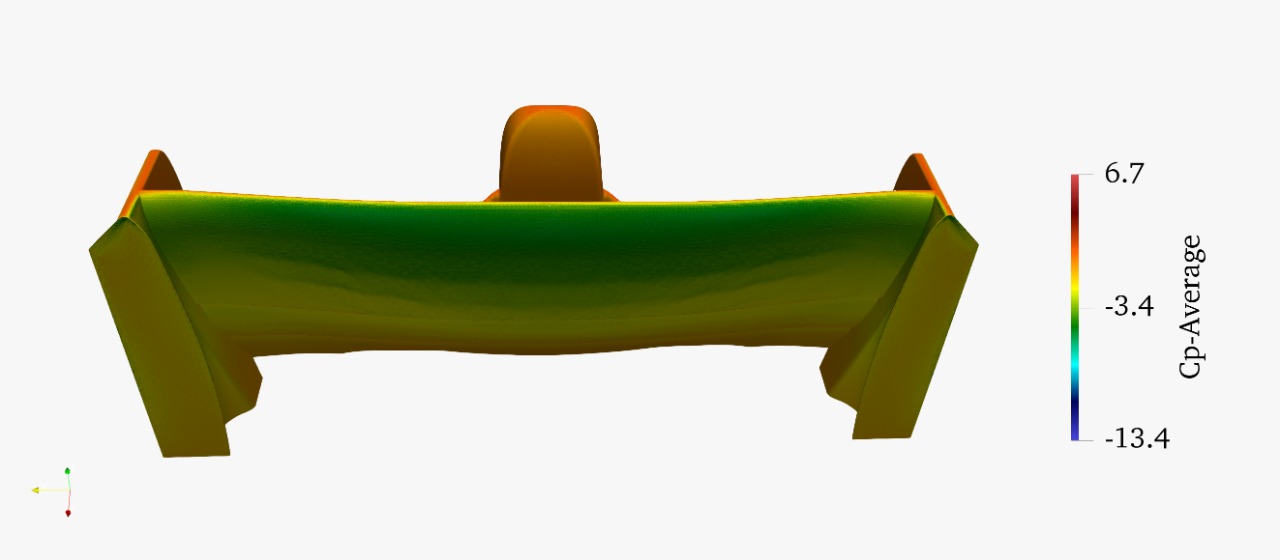}
    \caption{Bottom.}
  \end{subfigure}%
  \caption{Time-averaged pressure coefficient, $C_p$, on the surface of the Imperial Front Wing.}
  \label{fig:cp_ifw}
\end{figure}

In Figure~\ref{fig:cp_ifw}, we present the contour plot of the time-averaged pressure coefficient on the surface of the front wing.
The statistics have been obtained by averaging the solution for approximately five flow-through time units.

\begin{table}[!htb]
\sisetup{output-exponent-marker=}
\sisetup{scientific-notation=fixed, fixed-exponent=0}
\centering\small
\caption{Performance of different methods for error and CFL-based step size
         controllers:
         Number of function evaluations (\#FE), rejected steps (\#R), and
         wall-clock time in seconds for the Imperial Front Wing using
         polynomials of degree $p = 2$.}
\label{tab:ifw-p2}
\setlength{\tabcolsep}{0.75ex}
\begin{tabular*}{\linewidth}{@{\extracolsep{\fill}}c *2c r@{\hskip 0.5ex}rr@{\hskip 0.5ex}r@{\hskip 1ex}cr@{\hskip 0.5ex}r@{\hskip 1ex}c@{}}
  \toprule
  Scheme & $\beta$ & $\tol$/$\cfl$ & \multicolumn{1}{c}{\#FE} & \multicolumn{1}{c}{\#R} & Wall-clock time (s) \\
  \midrule

  \RK[BS]{3}{3}[][FSAL]          & $(0.60, -0.20, 0.00)$
     &  $\tol = 10^{-8}$ & $   3973$ & $    0$ & \tablenum{2.6999e+03} \\
   & &  $\cfl =  1.0   $ & $   7678$ &        & \tablenum{5.1884e+03} \smallskip\\
  \ssp43                         & $(0.55, -0.27, 0.05)$
     &  $\tol = 10^{-8}$ & $   2980$ & $    0$ & \tablenum{1.9715e+03} \\
   & &  $\cfl =  1.0   $ & $   4992$ &        & \tablenum{3.2989e+03} \smallskip\\
  \RK{3}{5}[\ESstarp][FSAL]      & $(0.70, -0.23, 0.00)$
     &  $\tol = 10^{-8}$ & $   3400$ & $    4$ & \tablenum{2.2778e+03} \\
   & &  $\cfl =  1.0   $ & $   6516$ &        & \tablenum{4.3623e+03} \smallskip\\
  \RK{4}{9}[\ESstarp]            & $(0.25, -0.12, 0.00)$
     &  $\tol = 10^{-8}$ & $   3195$ & $    1$ & \tablenum{2.1152e+03} \\
   & &  $\cfl =  1.0   $ & $   6111$ &        & \tablenum{4.0461e+03} \smallskip\\
  \RK{5}{10}[\ESstarp][FSAL]     & $(0.45, -0.13, 0.00)$
     &  $\tol = 10^{-8}$ & $   4094$ & $    3$ & \tablenum{2.7141e+03} \\
   & &  $\cfl =  1.0   $ & $   7811$ &         & \tablenum{5.1868e+03} \\

  \bottomrule
\end{tabular*}
\end{table}

The performance characteristics of the different methods for $p = 2$ are summarized
in Table~\ref{tab:ifw-p2}. The results are in agreement with those obtained
for the NASA juncture flow.
The CFL adaptor with $\cfl = 1.0$  works for all methods and is less efficient
than error-based step size controllers.
Again, \RK[BS]{3}{3}[][FSAL] is an efficient general purpose scheme for this
problem. Nevertheless, the optimized third- and fourth-order methods are more
efficient and \ssp43 is the most efficient scheme for this problem.
The fifth-order method is less efficient than the other four schemes, as expected.

\begin{table}[!htb]
\sisetup{output-exponent-marker=}
\sisetup{scientific-notation=fixed, fixed-exponent=0}
\centering\small
\caption{Performance of different methods for error and CFL-based step size
         controllers:
         Number of function evaluations (\#FE), rejected steps (\#R), and
         wall-clock time in seconds for the Imperial Front Wing using
         polynomials of degree $p = 1$. \ssp43 crashed for $\cfl = 1.0$.}
\label{tab:ifw-p1}
\setlength{\tabcolsep}{0.75ex}
\begin{tabular*}{\linewidth}{@{\extracolsep{\fill}}c *2c r@{\hskip 0.5ex}rr@{\hskip 0.5ex}r@{\hskip 1ex}cr@{\hskip 0.5ex}r@{\hskip 1ex}c@{}}
  \toprule
  Scheme & $\beta$ & $\tol$/$\cfl$ & \multicolumn{1}{c}{\#FE} & \multicolumn{1}{c}{\#R} & Wall-clock time (s) \\
  \midrule

  \RK[BS]{3}{3}[][FSAL]          & $(0.60, -0.20, 0.00)$
     &  $\tol = 10^{-8}$ & $   1411$ & $    0$ & \tablenum{3.5675e+02} \\
   & &  $\cfl =  1.0   $ & $   1579$ &        & \tablenum{3.9855e+02} \smallskip\\
  \ssp43                         & $(0.55, -0.27, 0.05)$
     &  $\tol = 10^{-8}$ & $   1352$ & $    0$ & \tablenum{3.3197e+02} \\
   & &  $\cfl =  0.85  $ & $   1304$ &        & \tablenum{3.2154e+02} \smallskip\\
  \RK{3}{5}[\ESstarp][FSAL]      & $(0.70, -0.23, 0.00)$
     &  $\tol = 10^{-8}$ & $   1127$ & $    6$ & \tablenum{2.8019e+02} \\
   & &  $\cfl =  1.0   $ & $   1341$ &        & \tablenum{3.3321e+02} \smallskip\\
  \RK{4}{9}[\ESstarp]            & $(0.25, -0.12, 0.00)$
     &  $\tol = 10^{-8}$ & $   1233$ & $    1$ & \tablenum{3.0231e+02} \\
   & &  $\cfl =  1.0   $ & $   1260$ &        & \tablenum{3.0834e+02} \smallskip\\
  \RK{5}{10}[\ESstarp][FSAL]     & $(0.45, -0.13, 0.00)$
     &  $\tol = 10^{-8}$ & $   1492$ & $    1$ & \tablenum{3.6616e+02} \\
   & &  $\cfl =  1.0   $ & $   1611$ &         & \tablenum{3.9535e+02} \\

  \bottomrule
\end{tabular*}
\end{table}

The results for $p = 1$ summarized in Table~\ref{tab:ifw-p1} are mostly
similar to the ones presented before. In contrast to the results for $p = 2$,
the CFL adaptor with $\cfl = 1.0$ did not work for \ssp43; the simulation
crashed when using $\cfl = 1.0$ and manual tuning was necessary to get a working
setup\footnote{This behavior can be explained by the different shape of the
stability region of \ssp43 compared to the other methods, with a relatively
larger real stability interval, \cf Table~1 in the supplementary material.
Another measure of the size of the stability region might help but would not remove
the necessity of manual tuning to get good performance for CFL-based controllers.}.
For $\cfl = 0.85$, the CFL adaptor worked and was a few percent more efficient
than the error-based controller. However, the latter did not require any manual
tuning at all and worked robustly with default parameters for all RK methods.
Here, \ssp43 is less efficient than \RK{3}{5}[\ESstarp][FSAL] and
\RK{4}{9}[\ESstarp]. Otherwise, the results are similar to the ones obtained
for the juncture flow and the setup using $p = 2$.

\section{Summary and conclusions}
\label{sec:summary}

We studied explicit Runge-Kutta methods applied to dissipative spectral element
semidiscretizations of hyperbolic conservation laws and CFD problems based on
the compressible Euler and Navier-Stokes equations.
In this context, we argued in Section~\ref{sec:cfl-vs-error} that error-based
step size control can be advantageous compared to CFL-based approaches, since
associated user-defined parameters are usually more robust and can be varied
in rather large ranges without affecting accuracy or efficiency.
Additionally, error-based step size control moves the burden of constructing
critical parts of the controller from the developer of the spatial
semidiscretization to the developer of the time integrator, easing the workflow
for most researchers.  The results for more complex test problems in Section
\ref{sec:numerical-experiments} also support this conclusion.

In Section~\ref{sec:importance-of-controller}, we demonstrated that choosing
good step size controller parameters is especially important if the time step
is restricted
by stability constraints, as is typical for many convection-dominated problems.
We compared existing Runge-Kutta pairs in
Section~\ref{sec:comparing-existing-schemes} and proposed an approach to
optimize controller parameters for such methods. In general, the third-order
method \RK[BS]{3}{3}[][FSAL] of Bogacki and Shampine \cite{bogacki1989a32}
performs well compared to both general purpose schemes and methods designed
specifically for CFD applications. The strong-stability preserving
method \ssp43 of \cite{kraaijevanger1991contractivity} with embedded method
of \cite{conde2018embedded} also performs well, but is more sensitive
to the choice of error tolerance.

In Section \ref{sec:optimized-RK}, we developed explicit low-storage Runge-Kutta pairs with optimized
step size controllers. These novel schemes are more efficient than all of the
existing schemes when applied to advection-dominated problems. We demonstrated
their performance in several CFD applications with increasing complexity,
including the compressible Euler and Navier-Stokes equations.
We contributed our optimized methods to the freely available open source
library DifferentialEquations.jl \cite{rackauckas2017differentialequations}
written in Julia \cite{bezanson2017julia}.

Although not demonstrated in this article,
another advantage of error-based step size control becomes apparent for a cold
startup of CFD problems, \ie simulations around complex geometries that are
initialized with a free stream flow. CFL-based approaches often need to adjust
the CFL scaling at the beginning to cope with the initial transient period.
In contrast, our error-based approach does not need special tuning and is robust
in our experience.

\begin{table}[tb]
\centering
\caption{Optimized PID controller parameters for some explicit Runge-Kutta
          methods with embedded error estimator.}
\label{tab:optimized=controllers}
\begin{tabular}{l *4c}
  \toprule
  Method & Reference & $\beta_1$ & $\beta_2$ & $\beta_3$ \\
  \midrule
  \RK[BS]{3}{3}[][FSAL] & \cite{bogacki1989a32}       & \num{0.60} & \num{-0.20} & \num{0.00} \\
  \RK[BS]{5}{7}[][FSAL] & \cite{bogacki1996efficient} & \num{0.28} & \num{-0.23} & \num{0.00} \\
  \RK[DP]{5}{6}[][FSAL] & \cite{prince1981high}       & \num{0.70} & \num{-0.40} & \num{0.00} \\
  \RK[T]{5}{6}[][FSAL]  & \cite{tsitouras2011runge}   & \num{0.57} & \num{-0.24} & \num{0.04} \\
  \ssp33  & \cite{shu1988efficient,conde2018embedded} & \num{0.70} & \num{-0.37} & \num{0.05} \\
  \ssp43  & \cite{kraaijevanger1991contractivity,conde2018embedded} & \num{0.55} & \num{-0.27} & \num{0.05} \\
  \RK[KCL]{3}{4}[\ZRp][][C] & \cite{kennedy2000low}    & \num{0.50} & \num{-0.35} & \num{0.10} \\
  \RK[KCL]{4}{5}[\ERp][][C] & \cite{kennedy2000low}    & \num{0.41} & \num{-0.28} & \num{0.08} \\
  \RK{3}{5}[\ESstarp] & this article          & \num{0.64} & \num{-0.31} & \num{0.04} \\
  \RK{3}{5}[\ESstarp][FSAL] & this article      & \num{0.70} & \num{-0.23} & \num{0.00} \\
  \RK{4}{9}[\ESstarp] & this article          & \num{0.25} & \num{-0.12} & \num{0.00} \\
  \RK{4}{9}[\ESstarp][FSAL] & this article      & \num{0.38} & \num{-0.18} & \num{0.01} \\
  \RK{5}{10}[\ESstarp] & this article         & \num{0.47} & \num{-0.20} & \num{0.06} \\
  \RK{5}{10}[\ESstarp][FSAL] & this article     & \num{0.45} & \num{-0.13} & \num{0.00} \\
  \bottomrule
\end{tabular}
\end{table}

A summary of existing and novel methods with optimized controller parameters
is given in Table~\ref{tab:optimized=controllers}.
Depending on whether dissipative/low-Mach effects dominate, \ssp43 and
\RK{3}{5}[\ESstarp][FSAL] are the most efficient schemes in our experience.
Additionally, \RK[BS]{3}{3}[][FSAL] is a surprisingly efficient general purpose
method.
It becomes increasingly complicated to design controllers that are stable
and efficient across different applications for methods of higher order and/or
with more stages. However, this is not necessarily a severe drawback, since
third-order accurate methods like \RK{3}{5}[\ESstarp][FSAL], \ssp43, and
\RK[BS]{3}{3}[][FSAL] are usually more efficient in CFD applications.
As argued in Section~\ref{sec:cfl-vs-error}, the error-based control has usually
a relatively mild sensitivity with respect to the choice of the tolerance.
In our experience, it is usually good to choose a relatively tight tolerance
around $10^{-8}$ for applied CFD problems. Since the step size is almost always limited by
stability, the tolerance does not matter that much, but a relatively tight
tolerance helps for methods that are more difficult to control (\eg those of
higher order or with more stages).

The present work was influenced and partially motivated by the landmark
work of Kennedy, Carpenter, and Lewis \cite{kennedy2000low}, which focused
on developing optimized Runge-Kutta methods for CFD.  Herein, we focus on modern
spectral element semidiscretizations that introduce dissipation at element interfaces, \eg using upwind
numerical fluxes. Hence, the stability regions of our methods focus also on
the negative real axis, whereas ``imaginary axis stability is a high priority
to the methods'' designed in \cite[p.~183]{kennedy2000low}.
Furthermore, we concentrate on the common case where the spatial error dominates
and the step size is restricted by stability rather than temporal accuracy.
Hence, we design Runge-Kutta pairs with large stability regions, for both the 
main and the embedded method. In particular, the stability regions of our novel
embedded schemes are larger than the ones of the corresponding main methods.

To our knowledge, this article is the first exploring the impact of controller
parameters and step size control stability on the efficiency of explicit
Runge-Kutta methods for CFD systematically. This provides important insights
into the construction of new methods and augments best practices published
before. In particular, we think the conventional wisdom that ``coping with step size
control instability is probably best accomplished by reducing step sizes''
\cite[p.~208]{kennedy2000low} can be improved upon by instead optimizing the
controller, since that results in a more robust and efficient scheme.
As noted in \cite[p.~208]{kennedy2000low}, ``doing this optimization requires
some caution because it is not sufficient in the design of a good controller
for each of the eigensolutions to be damped. The time constants associated
with these eigensolutions must not be too large or too small.''. Herein, we
proposed a way to conduct this optimization systematically and applied it
to a wide range of schemes.

Of course, such an approach also comes with limitations, in particular if
the main method is fixed such that only the embedded method and the controller
can be designed freely. Some methods such as the fourth-order method used in
this article constrain the range of embedded methods and controllers such that
a good general purpose optimization is not necessarily successful. While the
combined method can be efficient for certain problems, it is not necessarily
similarly efficient for other problems, \eg when going from inviscid to viscous
flows. Other schemes such as the novel third- and fifth-order accurate
optimized methods result in less stability restrictions, making the resulting
methods and controllers efficient for a broad range of CFD problems.
Thus, we would like to stress that designing a good time integration method
should not only focus on the main method but consider the interaction of
a main method, an error estimator, and a step size controller.
Applying this principle to viscous flows will be a subject of future work.

Some previous work has focused on automated step size control for
convection-dominated problems with the goal of achieving a temporal
error that is of similar magnitude to that of the spatial error
\cite{berzins1995temporal,ware1995adaptive}.  The addition of such a control
on top of the techniques employed here might lead to an even more efficient
controller that is not adversely affected by excessively tight temporal error
tolerance specification.

We expect the new methods developed in this article to perform well for
convection-dominated flows in the subsonic regime; here, we tested them
mainly with reference Mach numbers in the range 0.1--0.5.
We have demonstrated their improved performance compared to some standard schemes
also in other regimes, including viscous flows and the transonic/supersonic regime.
For small Mach numbers, incompressible solvers with implicit time discretizations
are usually applied. However, if compressible solvers should be used for low
Mach numbers, methods could be optimized following the
approach of this article.

\appendix
\section*{Acknowledgments}

Research reported in this publication was supported by the
King Abdullah University of Science and Technology (KAUST).
We are thankful for the computing resources of the
Supercomputing Laboratory and the Extreme Computing Research Center
at KAUST.
Funded by the Deutsche Forschungsgemeinschaft (DFG, German Research Foundation)
under Germany's Excellence Strategy EXC 2044-390685587, Mathematics Münster:
Dynamics-Geometry-Structure. Special thanks are extended to the McLaren F1 racing Team for providing data, CAD 
geometries, and setup of the Imperial Front Wing test case.

\section{Efficient implementation of \texorpdfstring{\ssp43}{SSP3(2)4}}
\label{sec:ssp43}

The Butcher coefficients of \ssp43 are
\begin{equation}
\label{eq:ssp43-butcher}
\renewcommand\arraystretch{1.2}
\begin{array}{c|cccc}
  0 &  &  &  & \\
  \nicefrac{1}{2} & \nicefrac{1}{2} &  &  & \\
  1 & \nicefrac{1}{2} & \nicefrac{1}{2} &  & \\
  \nicefrac{1}{2} & \nicefrac{1}{6} & \nicefrac{1}{6} & \nicefrac{1}{6} & \\
  \hline
  & \nicefrac{1}{6} & \nicefrac{1}{6} & \nicefrac{1}{6} & \nicefrac{1}{2}\\
  & \nicefrac{1}{4} & \nicefrac{1}{4} & \nicefrac{1}{4} & \nicefrac{1}{4}
\end{array}
\end{equation}
where spaces indicate zeros.
Because of its low-storage structure, the method can be implemented efficiently
and memory-friendly as
\begin{subequations}
\label{eq:ssp43-step}
\begin{equation}
  u
  \leftarrow
  u^n + \frac{1}{2} \dt_n f(t_n, u^n),
  \quad
  u
  \leftarrow
  u + \frac{1}{2} \dt_n f(t_n + \dt_n / 2, u),
  \quad
  u
  \leftarrow
  u + \frac{1}{2} \dt_n f(t_n + \dt_n, u),
\end{equation}
\begin{equation}
  \uhat
  \leftarrow
  \frac{1}{3} u^n + \frac{2}{3} u,
  \quad
  u
  \leftarrow
  \frac{2}{3} u^n + \frac{1}{3} u,
  \quad
  u
  \leftarrow
  u + \frac{1}{2} \dt_n f(t_n + \dt_n / 2, u),
  \quad
  \uhat
  \leftarrow
  \frac{1}{2} (\uhat + u).
\end{equation}
\end{subequations}
At the end of one step \eqref{eq:ssp43-step}, $u^{n+1}$ is stored in $u$
and $\uhat^{n+1}$ is stored in $\uhat$. Usually, it is not important to know
$\uhat^{n+1}$ but $\uhat^{n+1} - u^{n+1}$ to estimate the error; this difference
can be obtained as $(\uhat - u) / 2$ instead of $(\uhat + u) / 2$ in the last
assignment in \eqref{eq:ssp43-step}.
If the low-storage assumption introduced in \cite{ketcheson2010runge} can be
applied, \ssp43 can be implemented using only three memory locations for
$u^n$, $u$, and $\uhat$. Otherwise, an additional storage location is necessary
to evaluate the right-hand side $f$. Note that the previous value $u^n$ is
already included in this count of memory locations.

\section{Coefficients of the novel Runge-Kutta pairs}
\label{sec:coefficients}

The low-storage coefficients of the novel methods are listed in double
precision in Tables~\ref{tab:3Sstarp35}--\ref{tab:3SstarpFSAL510}.
Full-precision results in electronic form are available in the
accompanying repository \cite{ranocha2021optimizedCoefficients}.
We contributed our optimized methods to the freely available open source
library DifferentialEquations.jl \cite{rackauckas2017differentialequations}
written in Julia \cite{bezanson2017julia}.

\begin{table}[hbp]
\caption{Coefficients of the optimized, explicit, low-storage Runge-Kutta method
         \RK{3}{5}[\ESstarp].}
\label{tab:3Sstarp35}
\centering
\begin{small}%
\begin{tabular*}{\linewidth}{@{\extracolsep{\fill}} *3c @{}}
  \toprule
  $\gamma_{1,i}$ & $\gamma_{2,i}$ & $\gamma_{3,i}$ \\
  \midrule
  \verb|+0.0000000000000000e+00| & \verb|+1.0000000000000000e+00| & \verb|+0.0000000000000000e+00| \\
  \verb|+2.5876690703520788e-01| & \verb|+5.5284187451021605e-01| & \verb|+0.0000000000000000e+00| \\
  \verb|-1.3243668739945030e-01| & \verb|+6.7318444003896738e-01| & \verb|+0.0000000000000000e+00| \\
  \verb|+5.0556012314603993e-02| & \verb|+2.8031038045076351e-01| & \verb|+2.7525858134466369e-01| \\
  \verb|+5.6705528079028777e-01| & \verb|+5.5215088735073936e-01| & \verb|-8.9505487092797853e-01| \\
  \toprule
  $\delta_i$ & $\beta_{i}$ & $\bhat_i$ \\
  \midrule
  \verb|+1.0000000000000000e+00| & \verb|+1.1479315633699007e-01| & \verb|+1.0463633713540937e-01| \\
  \verb|+3.4076872093214550e-01| & \verb|+8.9335592952328596e-02| & \verb|+9.5204315749567586e-02| \\
  \verb|+3.4143992805846252e-01| & \verb|+4.3558587173792318e-01| & \verb|+4.4824466455686685e-01| \\
  \verb|+7.2293027328755899e-01| & \verb|+2.4735852952572862e-01| & \verb|+2.4490302954613102e-01| \\
  \verb|+0.0000000000000000e+00| & \verb|+1.1292684944702953e-01| & \verb|+1.0701165301202518e-01| \\
  \bottomrule
\end{tabular*}
\end{small}%
\end{table}

\begin{table}[hbp]
\caption{Coefficients of the optimized, explicit, low-storage Runge-Kutta method
         \RK{3}{5}[\ESstarp][FSAL].}
\label{tab:3SstarpFSAL35}
\centering
\begin{small}%
\begin{tabular*}{\linewidth}{@{\extracolsep{\fill}} *3c @{}}
  \toprule
  $\gamma_{1,i}$ & $\gamma_{2,i}$ & $\gamma_{3,i}$ \\
  \midrule
  \verb|+0.0000000000000000e+00| & \verb|+1.0000000000000000e+00| & \verb|+0.0000000000000000e+00| \\
  \verb|+2.5877719797257331e-01| & \verb|+5.5283549093013895e-01| & \verb|+0.0000000000000000e+00| \\
  \verb|-1.3243803601407234e-01| & \verb|+6.7318716082030616e-01| & \verb|+0.0000000000000000e+00| \\
  \verb|+5.0560339481908259e-02| & \verb|+2.8031039632976723e-01| & \verb|+2.7525632733046762e-01| \\
  \verb|+5.6705320007393134e-01| & \verb|+5.5215254470206099e-01| & \verb|-8.9505261746740339e-01| \\
  \toprule
  $\delta_i$ & $\beta_{i}$ & $\bhat_i$ \\
  \midrule
  \verb|+1.0000000000000000e+00| & \verb|+1.1479359710235412e-01| & \verb|+9.4841667050357029e-02| \\
  \verb|+3.4076558793345252e-01| & \verb|+8.9334428531133159e-02| & \verb|+1.7263713394303537e-01| \\
  \verb|+3.4143826550033862e-01| & \verb|+4.3558710250086169e-01| & \verb|+3.9982431890843712e-01| \\
  \verb|+7.2292753667879872e-01| & \verb|+2.4735761882014512e-01| & \verb|+1.7180168075801786e-01| \\
  \verb|+0.0000000000000000e+00| & \verb|+1.1292725304550591e-01| & \verb|+5.8819144221557401e-02| \\
  \bottomrule
\end{tabular*}
\end{small}%
\end{table}

\begin{table}
\caption{Coefficients of the optimized, explicit, low-storage Runge-Kutta method
         \RK{4}{9}[\ESstarp].}
\label{tab:3Sstarp49}
\centering
\begin{small}%
\begin{tabular*}{\linewidth}{@{\extracolsep{\fill}} *3c @{}}
  \toprule
  $\gamma_{1,i}$ & $\gamma_{2,i}$ & $\gamma_{3,i}$ \\
  \midrule
  \verb|+0.0000000000000000e+00| & \verb|+1.0000000000000000e+00| & \verb|+0.0000000000000000e+00| \\
  \verb|-4.6556413012591804e+00| & \verb|+2.4992627526078262e+00| & \verb|+0.0000000000000000e+00| \\
  \verb|-7.7202649248360644e-01| & \verb|+5.8668203654361373e-01| & \verb|+0.0000000000000000e+00| \\
  \verb|-4.0244232134197242e+00| & \verb|+1.2051413654126708e+00| & \verb|+7.6210371111381703e-01| \\
  \verb|-2.1296852467390187e-02| & \verb|+3.4747937967008691e-01| & \verb|-1.9811821590872183e-01| \\
  \verb|-2.4350225192344701e+00| & \verb|+1.3213461401287232e+00| & \verb|-6.2289607063175667e-01| \\
  \verb|+1.9856274809861678e-02| & \verb|+3.1196363243793707e-01| & \verb|-3.7522469934326264e-01| \\
  \verb|-2.8107901128852841e-01| & \verb|+4.3514190558940874e-01| & \verb|-3.3554365390009466e-01| \\
  \verb|+1.6894348958355357e-01| & \verb|+2.3596982994407883e-01| & \verb|-4.5609631107174843e-02| \\
  \toprule
  $\delta_i$ & $\beta_{i}$ & $\bhat_i$ \\
  \midrule
  \verb|+1.0000000000000000e+00| & \verb|+4.5037319691658841e-02| & \verb|+4.5506559279709452e-02| \\
  \verb|+1.2629238543878065e+00| & \verb|+1.8592173220119687e-01| & \verb|+1.1759683104926386e-01| \\
  \verb|+7.5749671775608729e-01| & \verb|+3.3297275092076306e-02| & \verb|+3.6582573305152133e-02| \\
  \verb|+5.1635911581112226e-01| & \verb|-4.7842226210501985e-03| & \verb|-5.3115558343556296e-03| \\
  \verb|-2.7463337920428273e-02| & \verb|+4.0558480626375678e-03| & \verb|+5.1782500127131271e-03| \\
  \verb|-4.3826746539417710e-01| & \verb|+4.1850279996827944e-01| & \verb|+4.9546390221186826e-01| \\
  \verb|+1.2735871036683928e+00| & \verb|-4.3818945074742778e-03| & \verb|-5.9993031327378659e-03| \\
  \verb|-6.2947400454427949e-01| & \verb|+2.7128460973244426e-02| & \verb|+9.4050934345683165e-02| \\
  \verb|+0.0000000000000000e+00| & \verb|+2.9522268113943101e-01| & \verb|+2.1693180876270352e-01| \\
  \bottomrule
\end{tabular*}
\end{small}%
\end{table}

\begin{table}
\caption{Coefficients of the optimized, explicit, low-storage Runge-Kutta method
         \RK{4}{9}[\ESstarp][FSAL].}
\label{tab:3SstarpFSAL49}
\centering
\begin{small}%
\begin{tabular*}{\linewidth}{@{\extracolsep{\fill}} *3c @{}}
  \toprule
  $\gamma_{1,i}$ & $\gamma_{2,i}$ & $\gamma_{3,i}$ \\
  \midrule
  \verb|+0.0000000000000000e+00| & \verb|+1.0000000000000000e+00| & \verb|+0.0000000000000000e+00| \\
  \verb|-4.6556414473350687e+00| & \verb|+2.4992627925744948e+00| & \verb|+0.0000000000000000e+00| \\
  \verb|-7.7202650996458722e-01| & \verb|+5.8668203777188754e-01| & \verb|+0.0000000000000000e+00| \\
  \verb|-4.0244366905198063e+00| & \verb|+1.2051460865230945e+00| & \verb|+7.6210066787213149e-01| \\
  \verb|-2.1296762840185311e-02| & \verb|+3.4747937221867325e-01| & \verb|-1.9811825043394005e-01| \\
  \verb|-2.4350225097901097e+00| & \verb|+1.3213460609651131e+00| & \verb|-6.2289592186990073e-01| \\
  \verb|+1.9856272971319869e-02| & \verb|+3.1196364646941938e-01| & \verb|-3.7522483807759566e-01| \\
  \verb|-2.8107911467910385e-01| & \verb|+4.3514195396843791e-01| & \verb|-3.3554383091351697e-01| \\
  \verb|+1.6894341687548597e-01| & \verb|+2.3596981300287537e-01| & \verb|-4.5609550050311212e-02| \\
  \toprule
  $\delta_i$ & $\beta_{i}$ & $\bhat_i$ \\
  \midrule
  \verb|+1.0000000000000000e+00| & \verb|+4.5037326272637540e-02| & \verb|+2.4836759124515911e-02| \\
  \verb|+1.2629238766481143e+00| & \verb|+1.8592173036998480e-01| & \verb|+1.8663277745621037e-01| \\
  \verb|+7.5749671896859117e-01| & \verb|+3.3297296725697173e-02| & \verb|+5.6710807959369842e-02| \\
  \verb|+5.1635894531407278e-01| & \verb|-4.7842041809589755e-03| & \verb|-3.4476954391492879e-03| \\
  \verb|-2.7463274218026097e-02| & \verb|+4.0558359610313108e-03| & \verb|+3.6022450565166364e-03| \\
  \verb|-4.3826731781279443e-01| & \verb|+4.1850277725960744e-01| & \verb|+4.5455706221450887e-01| \\
  \verb|+1.2735872946026565e+00| & \verb|-4.3819019689193264e-03| & \verb|-2.4346652894276124e-04| \\
  \verb|-6.2947402839274003e-01| & \verb|+2.7128437964460898e-02| & \verb|+6.6427553611035500e-02| \\
  \verb|+0.0000000000000000e+00| & \verb|+2.9522270159645919e-01| & \verb|+1.6136970795235051e-01| \\
  \bottomrule
\end{tabular*}
\end{small}%
\end{table}

\begin{table}
\caption{Coefficients of the optimized, explicit, low-storage Runge-Kutta method
         \RK{5}{10}[\ESstarp].}
\label{tab:3Sstarp510}
\centering
\begin{small}%
\begin{tabular*}{\linewidth}{@{\extracolsep{\fill}} *3c @{}}
  \toprule
  $\gamma_{1,i}$ & $\gamma_{2,i}$ & $\gamma_{3,i}$ \\
  \midrule
  \verb|+0.0000000000000000e+00| & \verb|+1.0000000000000000e+00| & \verb|+0.0000000000000000e+00| \\
  \verb|+4.0436600785046961e-01| & \verb|+6.8714670697523461e-01| & \verb|+0.0000000000000000e+00| \\
  \verb|-8.5034274642631846e-01| & \verb|+1.0930247604688987e+00| & \verb|+0.0000000000000000e+00| \\
  \verb|-6.9508941670724198e+00| & \verb|+3.2259753823301613e+00| & \verb|-2.3934051593421395e+00| \\
  \verb|+9.2387652253282782e-01| & \verb|+1.0411537008413965e+00| & \verb|-1.9028544220959867e+00| \\
  \verb|-2.5631780399574042e+00| & \verb|+1.2928214888647027e+00| & \verb|-2.8200422105832073e+00| \\
  \verb|+2.5457448699663476e-01| & \verb|+7.3914627692970059e-01| & \verb|-1.8326984641305650e+00| \\
  \verb|+3.1258317338631691e-01| & \verb|+1.2391292570393000e-01| & \verb|-2.1990945107506979e-01| \\
  \verb|-7.0071148005675854e-01| & \verb|+1.8427534793667669e-01| & \verb|-4.0824306603848765e-01| \\
  \verb|+4.8396209709807264e-01| & \verb|+5.7127889426970779e-02| & \verb|-1.3776697911212080e-01| \\
  \toprule
  $\delta_i$ & $\beta_{i}$ & $\bhat_i$ \\
  \midrule
  \verb|+1.0000000000000000e+00| & \verb|-2.2801023055963646e-03| & \verb|+5.7345884846761938e-02| \\
  \verb|-1.3317784091338497e-01| & \verb|+1.4073930208232305e-02| & \verb|+1.9714475180397338e-02| \\
  \verb|+8.2604227852460299e-01| & \verb|+2.3326917941728226e-01| & \verb|+7.2152966056837173e-02| \\
  \verb|+1.5137004305133324e+00| & \verb|+4.8082667004651816e-02| & \verb|+1.7396594898079398e-01| \\
  \verb|-1.3058100631770482e+00| & \verb|+4.1190032211396227e-01| & \verb|+3.7036936004454879e-01| \\
  \verb|+3.0366787893425076e+00| & \verb|-1.2914610713647529e-01| & \verb|-1.2155990390550650e-01| \\
  \verb|-1.4494582670745926e+00| & \verb|+1.2207460110385798e-01| & \verb|+1.1803729454911216e-01| \\
  \verb|+3.8343138733209576e+00| & \verb|+4.3578588031133875e-02| & \verb|+4.1556888233648698e-02| \\
  \verb|+4.1222939719233249e+00| & \verb|+1.0250768752899050e-01| & \verb|+1.2278866279103799e-01| \\
  \verb|+0.0000000000000000e+00| & \verb|+1.5593923403396062e-01| & \verb|+1.4562842322236844e-01| \\
  \bottomrule
\end{tabular*}
\end{small}%
\end{table}

\begin{table}
\caption{Coefficients of the optimized, explicit, low-storage Runge-Kutta method
         \RK{5}{10}[\ESstarp][FSAL].}
\label{tab:3SstarpFSAL510}
\centering
\begin{small}%
\begin{tabular*}{\linewidth}{@{\extracolsep{\fill}} *3c @{}}
  \toprule
  $\gamma_{1,i}$ & $\gamma_{2,i}$ & $\gamma_{3,i}$ \\
  \midrule
  \verb|+0.0000000000000000e+00| & \verb|+1.0000000000000000e+00| & \verb|+0.0000000000000000e+00| \\

  \verb|+4.0436601216857498e-01| & \verb|+6.8714670281614165e-01| & \verb|+0.0000000000000000e+00| \\
  \verb|-8.5034272895758400e-01| & \verb|+1.0930247489147509e+00| & \verb|+0.0000000000000000e+00| \\
  \verb|-6.9508941752621176e+00| & \verb|+3.2259753796071928e+00| & \verb|-2.3934051332441948e+00| \\
  \verb|+9.2387651927310854e-01| & \verb|+1.0411537025101014e+00| & \verb|-1.9028544224217609e+00| \\
  \verb|-2.5631780565098912e+00| & \verb|+1.2928214879121649e+00| & \verb|-2.8200422073999771e+00| \\
  \verb|+2.5457448793652260e-01| & \verb|+7.3914627557881230e-01| & \verb|-1.8326984652773810e+00| \\
  \verb|+3.1258317074119985e-01| & \verb|+1.2391292513718004e-01| & \verb|-2.1990944830846712e-01| \\
  \verb|-7.0071144144405084e-01| & \verb|+1.8427534723701233e-01| & \verb|-4.0824306358478707e-01| \\
  \verb|+4.8396210160238334e-01| & \verb|+5.7127889987965835e-02| & \verb|-1.3776697978802896e-01| \\
  \toprule
  $\delta_i$ & $\beta_{i}$ & $\bhat_i$ \\
  \midrule
  \verb|+1.0000000000000000e+00| & \verb|-2.2801003218369809e-03| & \verb|-2.0192554400120660e-02| \\
  \verb|-1.3317784195088034e-01| & \verb|+1.4073931157901863e-02| & \verb|+2.7379034809591845e-02| \\
  \verb|+8.2604228147502079e-01| & \verb|+2.3326917755084567e-01| & \verb|+3.0288186361459657e-01| \\
  \verb|+1.5137004257557283e+00| & \verb|+4.8082667413538623e-02| & \verb|-3.6568438806222223e-02| \\
  \verb|-1.3058100599350237e+00| & \verb|+4.1190032177069519e-01| & \verb|+3.9826647746767679e-01| \\
  \verb|+3.0366788029241634e+00| & \verb|-1.2914610678077362e-01| & \verb|-5.7159594211406851e-02| \\
  \verb|-1.4494582743988951e+00| & \verb|+1.2207460138487101e-01| & \verb|+9.8498551038485579e-02| \\
  \verb|+3.8343138991763621e+00| & \verb|+4.3578585831744204e-02| & \verb|+6.6546015524560853e-02| \\
  \verb|+4.1222937600129850e+00| & \verb|+1.0250768775680807e-01| & \verb|+9.0734795427481127e-02| \\
  \verb|+0.0000000000000000e+00| & \verb|+1.5593923423620598e-01| & \verb|+8.4322893253308037e-02| \\
  \bottomrule
\end{tabular*}
\end{small}%
\end{table}

\printbibliography

\includepdf[pages=-]{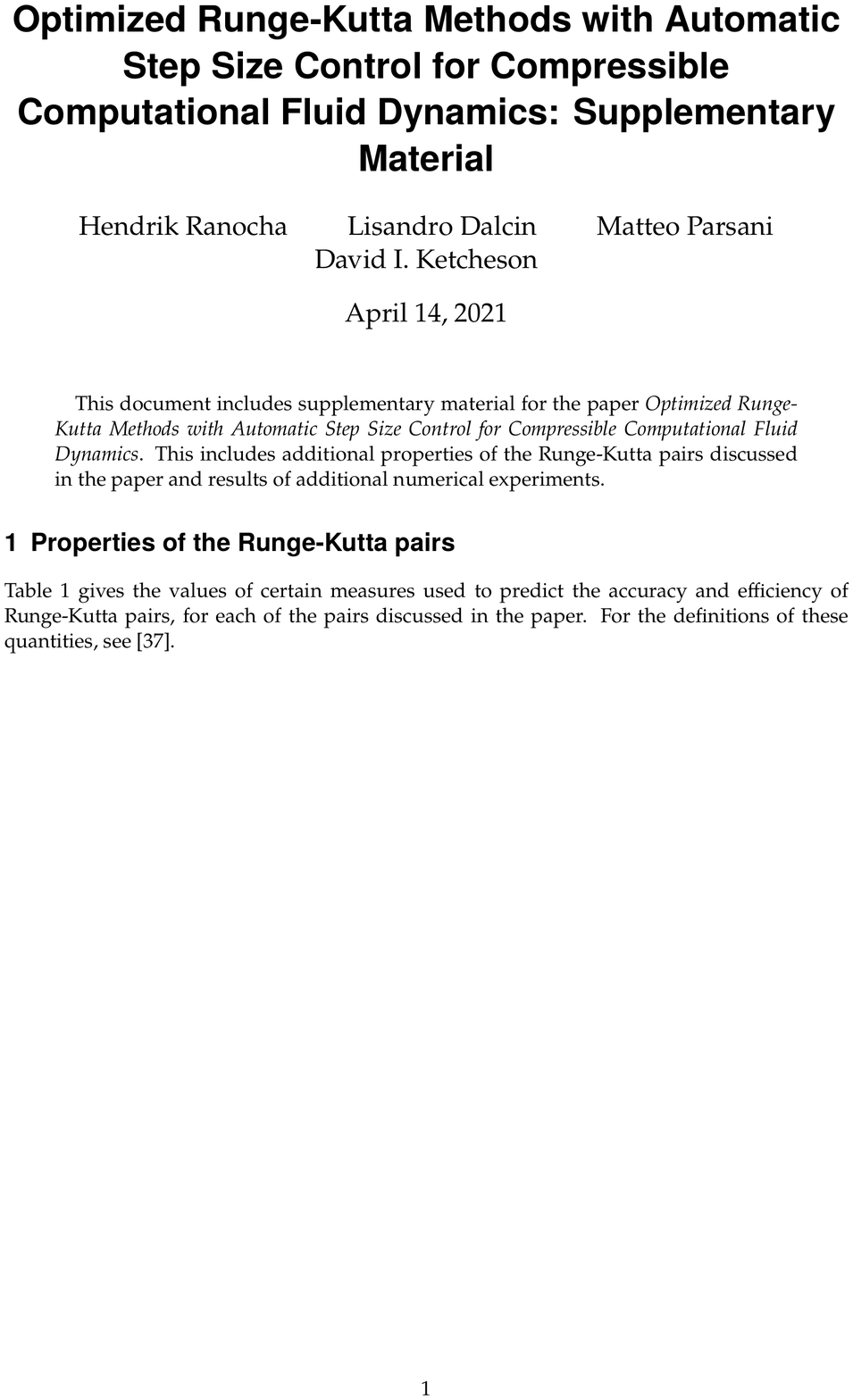}

\end{document}